%% file: DSBA.tex
\documentclass{article}

\usepackage{hyperref}       
\usepackage{url}            
\usepackage{booktabs}       
\usepackage{amsfonts}       
\usepackage{nicefrac}       
\usepackage{microtype}      
\usepackage{pifont}
\usepackage{amssymb}

\usepackage{algorithm}
\usepackage{algorithmic}
\usepackage{hyperref}
\usepackage{natbib}

\usepackage[a4paper,top=3cm,bottom=2cm,left=3cm,right=3cm,marginparwidth=1.75cm]{geometry}

\usepackage{amsmath,amssymb,amsthm}
\usepackage{pgfplots}
\usepackage{color}
\usepackage{flushend}
\usepackage{multirow}
\usepackage{rotating}
\usepackage{appendix}
\usepackage{tikz,pgfplots}
\usepackage{subcaption}
\input{mysymbol.sty}

\makeatletter

\newcommand{\IB}{\mathbf{I}}

\newcommand{\MB}{\mathbf{M}}

\newcommand{\UB}{\mathbf{U}}

\newcommand{\WB}{\mathbf{W}}

\newcommand{\zeroB}{\mathbf{0}}
\newcommand{\oneB}{\mathbf{1}}

\newcommand{\qB}{\mathbf{q}}

\newcommand{\zB}{\mathbf{z}}
\newcommand{\ZEROBB}{\mathbf{0}}

\newcommand{\BM}{\mathcal{B}}

\newcommand{\FM}{\mathcal{F}}

\newcommand{\TM}{\mathcal{T}}

\newcommand{\XM}{\mathcal{X}}
\newcommand{\YM}{\mathcal{Y}}

\newcommand{\NM}{\mathcal{N}}

\newcommand{\RM}{\mathcal{R}}

\newcommand{\RBB}{\mathbb{R}}

\newcommand{\find}{\mathop{\rm find}}

\newcommand{\defi}{\stackrel{\mathrm{def}}{=}}

\makeatother

\newtheorem{theorem}{Theorem}

\newtheorem{lemma}{Lemma}
\newtheorem{definition}{Definition}
\newtheorem{assumption}{Assumption}

\begin{document}

\title{A Decentralized Proximal Point-type Method for Saddle Point Problems}

\date{}


\author{
\begin{tabular}[t]{c@{\extracolsep{2em}}c @{\extracolsep{2em}}c} 
Weijie Liu\footnote{Authors are arranged in alphabetical order} & Aryan Mokhtari$^*$ & Asuman Ozdaglar$^*$\\
Zhejiang University & UT Austin & MIT \\
westonhunter@zju.edu.cn & mokhtari@austin.utexas.edu & asuman@mit.edu \\ \\ 
Sarath Pattathil$^*$ & Zebang Shen$^*$ & Nenggan Zheng$^*$ \\
MIT & University of Pennsylvania & Zhejiang University\\
sarathp@mit.edu & zebang@seas.upenn.edu & zng@cs.zju.edu.cn
\end{tabular}\\
\\
}

\maketitle

\begin{abstract}
In this paper, we focus on solving a class of constrained non-convex non-concave saddle point problems in a decentralized manner by a group of nodes in a network. Specifically, we assume that each node has access to a summand of a global objective function and nodes are allowed to exchange information only with their neighboring nodes. We propose a decentralized variant of the proximal point method for solving this problem. We show that when the objective function is $\rho$-weakly convex-weakly concave the iterates converge to approximate stationarity with a rate of $\mathcal{O}(1/\sqrt{T})$ where the approximation error depends linearly on $\sqrt{\rho}$. We further show that when the objective function satisfies the Minty VI condition (which generalizes the convex-concave case) we obtain convergence to stationarity with a rate of $\mathcal{O}(1/\sqrt{T})$. To the best of our knowledge, our proposed method is the first decentralized algorithm with theoretical guarantees for solving a non-convex non-concave decentralized saddle point problem. Our numerical results for training a general adversarial network (GAN) in a decentralized manner match our theoretical guarantees.
\end{abstract}

\input{1-introduction}

\input{2-preliminaries}

\input{3-problem}

\input{4-algorithm}

\input{5-theory}

\input{7-experiments}

\vspace{-1mm}
\section{Conclusion}
\label{sec:conclusion}
\vspace{-1mm}
In this paper we proposed DPPSP, the first distributed algorithm thats finds approximate first-order stationary points of a non-smooth non-convex non-concave saddle point problem. We showed that under the assumption that the objective function is $\rho$-weakly convex-weakly concave after $T$ iterations the iterates generated by DPPSP find a point that satisfy first order optimality up to an error of $\mathcal{O}(\sqrt{\rho}+1/\sqrt{T})$ while the gap between the iterates is of $\mathcal{O}(\sqrt{\rho}+1/\sqrt{T})$. We further improved the our  first order optimality complexity to $\mathcal{O}(1/\sqrt{T})$ and our consensus violation error to $\mathcal{O}(1/\sqrt{T})$ under the Minty Variational Inequality assumption. We also implemented DPPSP for training a general adversarial network and illustrated its superior performance in practice. 

\input{6-appendix}

\bibliographystyle{plainnat}
\bibliography{DSBA_minmax,bmc_article,bmc_article2,ex_bib,references}

\end{document}

%% file: 1-introduction.tex

\section{Introduction}
\label{sec:intro}

In this paper we focus on solving the following saddle point problem
\begin{align}
\min_{\bbx \in \mathcal{X}} \max_{\bby \in \mathcal{Y}} f(\bbx, \bby)
\label{prob:main_not_dist}
\end{align}
where $\bbx \in \mathcal{X}\subset \reals^p$ and $\bby \in \mathcal{Y}\subset \reals^q$ are the inputs of the function $f:\reals^{p\times q}\to \reals$. We assume that both sets $\mathcal{X}$ and $\mathcal{Y}$ are nonempty convex and compact. The min-max optimization problem in~\eqref{prob:main_not_dist} captures a wide range of important problems in game theory \cite{basar1999dynamic}, control theory \cite{hast2013pid}, robust optimization \cite{ben2009robust}, and more recently machine learning such as generative adversarial networks (GAN) \cite{goodfellow2014generative}, adversarial training \cite{madry2017towards}, and multi-agent reinforcement
learning \cite{omidshafiei2017deep}.

In this paper we are interested in studying the general non-convex non-concave version of Problem~\eqref{prob:main_not_dist} in which $f$ could be non-convex with respect to the minimization variable $\bbx$ and non-concave with respect to the maximization variable $\bby$. Indeed, in such case, finding an optimal solution of Problem~\eqref{prob:main_not_dist} is hard, in general, as we know solving a special case of this problem which is minimizing a non-convex function is hard, in general. Therefore, we settle for finding a set of points $(\bbx^*,\bby^*)$ that satisfy the first-order optimality condition for Problem~\eqref{prob:main_not_dist}, i.e., 
\begin{align}\label{eq_goal_global}
    \left( \nabla_\bbx f(\bbx^*, \bby^*)\right)^\top (\bbx-\bbx^*) &\geq 0 , \nonumber \\
        \left(\nabla_\bby f(\bbx^*, \bby^*)\right)^\top  (\bby-\bby^*) &\leq 0,
\end{align}
for any $\bbx\in \ccalX$ and $\bby\in \ccalY$. The focus of this paper is on distributed optimization in which the objective function $f$ is defined as a sum of $N$ individual functions $f_{n}:\reals^{p\times q}\to \reals$ for $n\in\{1,\dots,N\}$, i.e., $f(\bbx, \bby):=\sum_{n=1}^N f_n(\bbx, \bby)$. We assume that each component function $f_n$ is assigned to a node (machine), indexed by $n$, in a network of size $N$. Indeed, the considered distributed optimization setting is of interest due to the limited capacity of processors that we have access to as well as the massive-scale datasets that we aim to train. For instance, if the function $f$ represents a loss function over training data which is stored in different machines, each $f_i$ would represent the loss on the data in machine $i$ and the function $f = \sum f_i$ represents the overall training loss. We would like to have an algorithm which updates parameters in such a manner that no machine will need access to data from all other machines to update the parameters. Our goal is to ensure that all nodes in the network find a solution $(\bbx^*,\bby^*)\in \ccalX\times \ccalY$ that (approximately) satisfies the condition in \eqref{eq_goal_global} without any global coordination and only by exchanging information with their neighboring nodes. This is a challenging task as we need to ensure that the local variables of each node satisfies the condition in \eqref{eq_goal_global}, and at the same time they are required to equal each other so that the consensus condition is satisfied.

In this paper, we formally propose a fully distributed algorithm that finds a set of local iterates $\{\bbx_1,\dots,\bbx_N\}$ and $\{\bby_1,\dots,\bby_N\}$ which (approximately) satisfy the first-order stationarity condition for the saddle point problem in \eqref{prob:main_not_dist} and also (approximately) satisfy the consensus condition $\{\bbx_1=\dots=\bbx_N\}$ and $\{\bby_1=\dots=\bby_N\}$. To the best of our knowledge, our proposed method is the first \textit{distributed} algorithm that finds a \textit{first-order stationary point} of the min-max problem in \eqref{prob:main_not_dist} for the considered \textit{nonsmooth constrained non-convex non-concave} setting in a finite number of iterations. The rest of the paper is organized as follows. We first state recap some required definitions and assumptions for introducing our method and its analysis (Section~\ref{sec:prelim}). We then proceed by formally defining the distributed consensus problem that we aim to solve and its equivalent problem in matrix form (Section~\ref{sec:prob_formulation}). The proposed Decentralized Proximal Point method for Saddle Point (DPPSP) problems is then presented in both global and local forms (Section~\ref{sec:algorithm}). Then, the convergence analysis of DPPSP is provided (Section~\ref{sec:results}). Finally, numerical experiments for training a generative adversarial network is presented in (Section~\ref{sec:exp}). We close the paper by concluding remarks (Section~\ref{sec:conclusion}). The proofs are provided in the supplementary material.

\subsection{Related Work}

\textbf{Saddle-point problem.} Min-Max optimization, also known as saddle point problem, is a well studied area and several algorithms have been proposed to solve this problem. The celebrated Proximal Point method was introduced in \cite{martinet1970breve} and analyzed in \cite{rockafellar1976monotone} for the case that $f(\bbx, \bby)$ is strongly convex-strongly concave or bilinear. Following this, several inexact versions of the proximal method were proposed to solve the problem including the hybrid inexact proximal method \cite{monteiro2010complexity}, which generalizes the extragradient method studied in \cite{korpelevich1976extragradient, nemirovski2004prox}. Recently optimistic gradient descent ascent (OGDA) was proposed and analyzed for solving saddle point problems \cite{liang2018interaction, gidel2018variational,mokhtari2019unified,mokhtari2019proximal}. All these works consider the setting where $f(\bbx, \bby)$ is (strongly)convex-concave. Several other papers study the setting where the objective function $f(\bbx, \bby)$ is nonconvex with respect to $\bbx$ but is concave with respect to $\bby$ \cite{rafique2018non,nouiehed2019solving,lu2019hybrid}, or study specific cases of non-convex non-concave problems \cite{sanjabi2018solving,lin2018solving}. In a parallel line of work, several papers including \cite{nemirovski2009robust, chen2014optimal,palaniappan2016stochastic} study the stochastic version of the problem where we only have access to an unbiased estimate of the gradient and not the true gradient itself. 
For the general non-convex non-concave setting, most of the existing works focus on showing that the gradient descent ascent (GDA) method and its infinitesimal counterpart (GDA-dynamic) converge asymptotically to the local Nash equilibrium points \cite{nagarajan2017gradient,daskalakis2018limit,adolphs2018local,jin2019minmax}. 
However, all the mentioned works aim at finding a saddle point in a centralized setting and cannot be applied to decentralized settings.

\textbf{Decentralized optimization.} There are several works which study the problem of decentralized minimization when the objective function is (strongly) convex. For such settings, methods like Decentralized Gradient Descent (DGD) \cite{Nedic2009,Jakovetic2014-1,yuan2016convergence}, Augmented Lagrangian Method (ALM) \cite{shi2015extra,shi2015proximal,mokhtari2016dsa,shen2018towards}, distributed implementations of the alternating direction method of multipliers (ADMM) \cite{Schizas2008-1,BoydEtalADMM11,Shi2014-ADMM,chang2015multi}, decentralized dual averaging \cite{Duchi2012,cTsianosEtal12}, and several dual based strategies \cite{seaman2017optimal,scaman2018optimal,uribe2018dual} are proposed and their corresponding convergence guarantees are established. Recently, there have been some works which look at the problem of decentralized minimization when the function is non-convex and show convergence to a first-order stationary point \cite{zeng2018nonconvex, Hong_Prox-PDA, hong2018gradient, sun2018distributed, scutari2017parallel, scutari2018distributed}.
 However, all these works rely crucially on the fact that the goal is minimizing a function, and the analysis cannot be easily extended to the setting of min-max optimization which we consider in this paper. 

\textbf{Decentralized saddle-point problem.} There are multiple works which look at the decentralized saddle point problem when the objective function is convex-concave \cite{wai2018multi, mateos2015distributed, koppel2015saddle}. However, none of these works provide any convergence guarantees for non-convex non-concave problems.

\medskip\noindent{\bf Notation.\quad} Lowercase boldface $\bbv$ denotes a vector and uppercase boldface $\bbA$ denotes a matrix. We use $\|\bbv\|$ to denote the Euclidean norm of vector $\bbv$. Given a multi-input function $f(\bbx,\bby)$, its gradient with respect to $\bbx$ and $\bby$ at points $(\bbx_0,\bby_0)$ are denoted by $\nabla_\bbx f(\bbx_0,\bby_0)$ and $\nabla_\bby f(\bbx_0,\bby_0)$, respectively. We refer to the largest and smallest eigenvalues of a matrix $\bbA$ by $\lambda_{\max}(\bbA)$ and  $\lambda_{\min}(\bbA)$, respectively. We use the notation $\bbA\otimes \bbB$ to indicate the Kronecker product of matrices $\bbA$ and $\bbB$. We use $\bbone_p$ to denote a vector of size $p$ where all its components are $1$ and use $\bb0_p$ to denote a vector of size $p$ where all its components are $0$. The identity matrix of size $d\times d$ is denoted by $\bbI_d$.

%% file: 2-preliminaries.tex

\section{Preliminaries}\label{sec:prelim}
In this section, we present definitions and assumptions that we will be using throughout the paper. 
\begin{definition}
Consider a function $\phi:\reals^p\to \reals$. The function $\phi$ is\\
(a) convex over the set $\ccalX$ if for any $\bbx,\hbx\in \ccalX$ we have
\begin{align*}
\phi(\hbx) \geq \phi(\bbx) + \langle \nabla \phi(\bbx), \hbx - \bbx \rangle.
\end{align*}
(b) $\mu$-strongly convex over the set $\ccalX$ if for any $\bbx,\hbx\in \ccalX$ we have
\begin{align*}
\phi(\hbx) \geq \phi(\bbx) + \langle \nabla \phi(\bbx), \hbx - \bbx \rangle + \frac{\mu}{2} \| \hbx - \bbx\|^2 .
\end{align*}
(c) $\rho$-weakly convex over the set $\ccalX$ if for any $\bbx,\hbx\in \ccalX$ we have if 
\begin{align*}
\phi(\hbx) \geq \phi(\bbx) + \langle \nabla \phi(\bbx), \hbx - \bbx \rangle - \frac{\rho}{2} \| \hbx - \bbx\|^2 .
\end{align*}
Further, a function $\phi(\bbx)$ is concave, $\mu$-strongly concave, or $\rho$-weakly concave, if $-\phi(\bbx)$ is convex, $\mu$-strongly convex, or $\rho$-weakly convex, respectively.
\end{definition}

Note that a $\rho$-weakly convex function is a non-convex function with a specific structure that by adding a quadratic regularization term of $(\rho/2)\|\bbx\|^2$ it becomes convex.



\begin{definition}
Consider an operator $F:\reals^{d}\to\reals^d$. The operator $F$ is\\
(a) monotone over the set $\ccalZ$ if for any $\bbz,\hbz\in \ccalZ$ we have
\begin{align}
\langle F(\bbz) - F(\hbz), \bbz - \hbz \rangle \geq 0
\end{align}
(b) $\mu$-strongly monotone over the set $\ccalZ$ if for any $\bbz,\hbz\in \ccalZ$ we have 
\begin{align}
\langle F(\bbz) - F(\hbz), \bbz - \hbz \rangle \geq \mu \| \bbz - \hbz \|^2 
\end{align}
(c) $\rho$-weakly monotone over the set $\ccalZ$ if for any $\bbz,\hbz\in \ccalZ$ we have
\begin{align}
\langle F(\bbz) - F(\hbz), \bbz- \hbz \rangle \geq -\rho \| \bbz - \hbz \|^2 
\end{align}
\end{definition}


\begin{definition}
We call the inverse operator $\ccalG^{-1}$ \textit{well-defined}, if for any $\bbz$ the problem of finding $\bbu$ such that $\bbu=\ccalG^{-1}(\bbz)$ has a unique solution. As a consequence, the operator $(\ccalI+\ccalF)^{-1}$ is well-defined if the norm of the operator $\ccalF$ is strictly smaller than $1$, i.e., $\|\ccalF(\bbz)\|< \|\bbz\|$ for any $\bbz$.
\end{definition}
Now we are at the right point to state the assumptions that we will use in the rest of the paper.

\begin{assumption}\label{ass:weak_con}
The objective function $f(\bbx, \bby)$ is $\rho$-weakly convex with respect to $\bbx$ and $\rho$-weakly concave with respect to $\bby$ for all $\bbx\in \ccalX$ and $\bby \in \ccalY$.
\end{assumption}


\begin{assumption}\label{ass:compact}
The sets $\mathcal{X}$ and $\mathcal{Y}$ are convex, closed and bounded.
\end{assumption}

In the following lemma, we characterized the properties of $F([\bbx,\bby]) = [\nabla_{\bbx} f(\bbx, \bby) ; -\nabla_{\bby} f(\bbx, \bby)]$ under the considered assumptions. 
\begin{lemma}
Consider the variable $\bbz = [\bbx; \bby]$ and the operator $F(\bbz) = [\nabla_{\bbx} f(\bbx, \bby) ; -\nabla_{\bby} f(\bbx, \bby)]$, and the set $\ccalZ=\ccalX \times \ccalY$. If the conditions in Assumptions~\ref{ass:weak_con}-\ref{ass:compact} are satisfied then the operator $F$ is $\rho$-weakly monotone over the set $\ccalZ$, and $\ccalZ$ is convex, closed and bounded with a diameter $D$, i.e., for $\|\bbz\|\leq D$ for any $\bbz\in \ccalZ$.
\end{lemma}

%% file: 3-problem.tex


\section{Problem Formulation}\label{sec:prob_formulation}
\vspace{-1mm}
In this section, we first introduce a formulation of the decentralized version of the general saddle point problem in \eqref{prob:main_not_dist}. Consider a network of $N$ nodes (machines) in which each node $n$ has access to a local objective function $f_n: \reals^{p\times q}\to \reals$. Further, let $\bbx_n,\bby_n$ be the decision variables for node $n$. Our goal is to solve the following distributed consensus saddle point problem 
\begin{align}\label{eqn_centralized_consensus_minmax_problem}
	\min_{\{\bbx_n\}_{n=1}^N} \  & \max_{\{\bby_n\}_{n=1}^N}\ \sum_{n=1}^N f_n(\bbx_n, \bby_n)\nonumber \\
	\st\ \bbx_1  = \ldots  = \bbx_N & \in  \ccalX , \quad \bby_1 = \ldots = \bby_N\in \ccalY,
\end{align}
where to ensure that nodes find the same stationary point and reach consensus we enforce the iterates of the nodes to be equal to each other, i.e., $\bbx_1 = \ldots = \bbx_N$ and $\bby_1 = \ldots = \bby_N$. It can be easily verified that the distributed problem in \eqref{eqn_centralized_consensus_minmax_problem} is equivalent to the original problem in \eqref{prob:main_not_dist} when we enforce the consensus condition. In other words, a pair $(\bbx^*,\bby^*)$ is an optimal solution of \eqref{prob:main_not_dist}  if and only if $\bbx_1 = \ldots = \bbx_N=\bbx^*$ and $\bby_1 = \ldots = \bby_N=\bby^*$ is an optimal solution of  \eqref{eqn_centralized_consensus_minmax_problem}. As we described in \eqref{eq_goal_global}, solving the problem exactly is hard in general and our goal is to ensure that all nodes in the network find a point $(\bbx^*, \bby^*)$ that satisfies the first-order stationarity, i.e., for any $\bbx\in\ccalX$ and $\bby\in \ccalY$ we have
\begin{align}\label{eq_goal}
    \left(\sum_{n=1}^N \nabla_\bbx f_n(\bbx^*, \bby^*)\right)^\top (\bbx-\bbx^*) &\geq 0 , \nonumber \\
        \left(\sum_{n=1}^N \nabla_\bby f_n(\bbx^*, \bby^*)\right)^\top  (\bby-\bby^*) &\leq 0.
\end{align}
Let us define $\oneB_{\XM}$ and $\oneB_{\YM}$ as the indicator functions of the feasible sets $\XM$ and $\YM$, respectively, and use the notation $\partial \oneB_{\XM}$ and $\partial \oneB_{\YM}$ to denote their corresponding sub-gradients.
Based on a standard convex optimization argument $\partial \oneB_{\XM}$ and $\partial \oneB_{\YM}$ are the normal cones of the convex sets $\XM$ and $\YM$, respectively. Therefore, we can express the conditions in \eqref{eq_goal} as
\begin{align}\label{eq_goal_equal}
0&\in \sum_{n=1}^N \nabla_\bbx f_n(\bbx^*, \bby^*) + \partial \oneB_{\XM}(\bbx^*) \nonumber \\
0&\in - \sum_{n=1}^N \nabla_\bby f_n(\bbx^*, \bby^*) + \partial \oneB_{\YM}(\bby^*) .
\end{align}
Now define the local operators $\ccalB_n: \reals^{d} \to \reals^{d}$ and $\ccalR_n: \reals^{d} \to \reals^{d}$ at node $n$ as $\ccalB_n([\bbx;\bby])=[\nabla_\bbx  f_n(\bbx,\bby), -\nabla_\bby  f_n(\bbx,\bby)] $ and $\RM_n([\bbx;\bby]) = [\partial \oneB_{\XM}, \partial \oneB_{\YM}]$, where to simplify our notation we defined $d:=p+q$. Considering these definitions our goal in \eqref{eq_goal} and \eqref{eq_goal_equal} can be written as the following decentralized root finding problem
\begin{align}\label{eqn_decentralized_consensus_root_monotone_operator}
	\find_{\{\zB_n\}_{n=1}^N}\quad \st  \quad & \sum_{n=1}^{N} \BM_n(\zB_n) + \RM_n(\zB_n) = 0, \nonumber \\
	  & \ {\zB_1 = \ldots = \zB_N}\in \ccalZ,
\end{align}
where $\zB_n \defi [\bbx_n; \bby_n]\in \reals^{d}$ is the concatenated local variable at node $n$ and the set $\ccalZ\subset \reals^{d}$ is a convex defined as $\ccalZ:=\{[\bbx; \bby]\in \reals^{d}\mid \bbx\in\ccalX, \bby\in\ccalY\}$. 
 To further simplify our problem formulation and write it in a more compact form, let us define the vector $\bbz=[\zB_1;\dots;\zB_n] \in \reals^{Nd}$ as the concatenation of all the local variables, ${\ccalB(\bbz):= [\BM_1(\zB_1); \cdots; \BM_N(\zB_N)] } \in \reals^{Nd}$ as the operator corresponding to the aggregate objective function, ${\ccalR(\bbz):= [\ccalR_1(\zB_1); \cdots ; \ccalR_N(\zB_N)] } \in \reals^{Nd}$ as the aggregate operator corresponding to the indicator functions, and $\ccalZ^N\subset \reals^{Nd}$ is given by $\ccalZ^N:=\{[\bbz_1;\dots; \bbz_N]\in \reals^{Nd}\mid \bbz_1,\dots,\bbz_n\in\ccalZ\}$. Then, the problem in \eqref{eqn_decentralized_consensus_root_monotone_operator} can be stated as 
\begin{align}\label{eqn_decentralized_consensus_root_monotone_operator_version_2}
\find_{\zB}\quad  \st\  & (\bbone_N\otimes \bbI_d)^\top(\BM(\zB) + \RM(\zB))= 0,  \nonumber \\
 &(\IB_{Nd}-\hbW )\zB=\zeroB_{Nd},\quad   \bbz\in \ccalZ^{N},
\end{align}
where we replaced the consensus constraint $\zB_1 = \ldots = \zB_N$ with the condition $(\bbI-\hbW)\bbz=\bb0$. Here the matrix $\hbW = \WB \otimes \IB_{d}\in \reals^{Nd\times Nd}$ is the Kronecker product of the identity matrix $\bbI_d\in \reals^{d\times d}$ and a mixing matrix $\bbW\in\reals^{N\times N}$, where $\bbW$ has the sparsity pattern of the graph and is designed in a way that the constraint $(\bbI-\hbW)\bbz=\bb0$ is satisfied if and only if $\zB_1 = \ldots = \zB_N$. This is a common approach in decentralized optimization \cite{yuan2016convergence}, and if $\bbW$ satisfies the following conditions
\begin{align}\label{w_conditions}
&\bbW=\bbW^T, \quad \bbW\oneB_N=\oneB_N, \nonumber \\
  \mathrm{null} (\IB  - \WB & ) = \mathrm{span} \{\oneB_N\}, \quad \bb0_N \prec \bbW\preceq  \bbI_N,
\end{align}
then \eqref{eqn_decentralized_consensus_root_monotone_operator} and \eqref{eqn_decentralized_consensus_root_monotone_operator_version_2} are equivalent. 
The equality conditions in \eqref{eqn_decentralized_consensus_root_monotone_operator_version_2} are coupled as $\bbz$ has to be chosen for both of them at the same time. We proceed to define a new variable $\bbq$ (which behaves as the dual variable for the consensus constraint) to separate these equality conditions. If we define $\UB\triangleq(\IB-\hbW)^{1/2}$, then the optimality conditions of Problem~\eqref{eqn_decentralized_consensus_root_monotone_operator_version_2} imply that there exists some $\bbp^* \in \RBB^{Nd}$, such that for $\bbq^* = \UB \bbp^*\in \RBB^{Nd}$ and $\alpha>0$ we have	
\begin{equation}
		\UB \bbq^* + \alpha [\BM(\bbz^*) + \RM(\bbz^*)] = \zeroB ~\mathrm{and}~ -\UB \bbz^* = \zeroB,
		\label{eqn_optimality}
	\end{equation}
where $\bbz^*\in \RBB^{Nd}$ is a solution of Problem~\eqref{eqn_decentralized_consensus_root_monotone_operator_version_2}; see, e.g., \cite{shen2018towards}. Equation \eqref{eqn_optimality} is obtained from the first order stationarity conditions of Problem \eqref{eqn_decentralized_consensus_root_monotone_operator_version_2} where $\frac{1}{\alpha}\bbp^*$ can be seen as the Lagrange multiplier for the constraint in \eqref{eqn_decentralized_consensus_root_monotone_operator_version_2}. The first equation follows from the optimality condition and the second one follows from the constraint, since $-\bbU \bbz^* = \bb0 \iff (\bbI - \hat{\bbW})\bbz^* = 0$. Hence, instead of solving \eqref{eqn_decentralized_consensus_root_monotone_operator_version_2} we find the optimal pair $(\bbz^*,\bbq^*)$ for \eqref{eqn_optimality} by updating $\bbz$ and $\bbq$ alternatively, as we do in the following section.

%% file: 4-algorithm.tex


\section{Algorithm}\label{sec:algorithm}
In this section, we aim to design a fully decentralized algorithm for solving the root finding problem in \eqref{eqn_optimality} which leads to a set of local iterates satisfying the conditions in \eqref{eq_goal_equal}. If we define $\bbv=[\bbz;\bbq]\in \RBB^{2Nd}$ as the concatenation of $\bbz$ and $\bbq$, then our problem is equivalent to finding the root of $\ccalT(\bbv)$, where the operator $\ccalT:\reals^{2Nd}\to\reals^{2Nd}$ is defined as
\begin{equation}
\TM(\bbv) = \bigg(
\underbrace{\begin{bmatrix}
	\alpha[\BM + \RM]& 0 \\
	0& 0
	\end{bmatrix}}_{\TM_1} +
\underbrace{\begin{bmatrix}
	\zeroB & \UB \\
	-\UB& \zeroB
	\end{bmatrix}}_{\TM_2}
\bigg)
\underbrace{\begin{bmatrix}
	\bbz \\
	\bbq
	\end{bmatrix}}_{\bbv},
\label{eqn_augmented_operator}
\end{equation}
i.e., finding a point $\bbv^*$ such that $\ccalT(\bbv^*)=\bb0$. Under the condition that the operator $(\ccalI+\ccalT)^{-1}$ is well-defined, finding a root of $\ccalT$ is equivalent to finding a fixed point of the operator $(\ccalI+\ccalT)^{-1}$, which is a point that satisfies $\bbv^* =(\ccalI+\ccalT)^{-1}(\bbv^*)$. This problem can be solved by following the recursive update $\bbv^{t+1} =  (\ccalI+\ccalT)^{-1} (\bbv^{t})$. However, it can be verified that implementation of this algorithm in a distributed setting is infeasible as computing the inverse operator $(\ccalI+\ccalT)^{-1}$ requires global communication. To better highlight this point, consider the simplified case that $\alpha=0$. The operator $(\ccalI+\ccalT)^{-1}$ in this case is $[\bbI,\bbU;-\bbU,\bbI]^{-1}$ and computation of this inverse requires access to the operator $(\bbI-\bbU^2)^{-1}=\hbW^{-1}$ which cannot be implemented in a distributed way. 

To resolve this issue, we introduce a system which has the same root as $\ccalT(\bbv)=\bb0$ and can be implemented in a distributed fashion (by exchanging information only with neighboring nodes). To do so, we consider the problem of finding a fixed point of the operator $(\bbD+\ccalT)^{-1}\bbD$ instead of $(\ccalI+\ccalT)^{-1}$, where $\bbD\succ \bb0$ is a positive definite matrix. Note that if $\bbv^*$ is a fixed point of $(\bbD+\ccalT)^{-1}\bbD$, i.e., $\bbv^*=(\bbD+\ccalT)^{-1}\bbD(\bbv^*)$ then it  satisfies the condition $\bbD\bbv^*+\ccalT\bbv^*=\bbD\bbv^*$ which implies that $\bbv^*$ is a root of $\ccalT$. Therefore, if the operator $(\bbD+\ccalT)^{-1}$ is well-defined, 
by updating the iterates based on the following fixed point iteration 
\begin{equation}
	\bbv^{t+1}=(\bbD+\ccalT)^{-1}\bbD\bbv^t
	\label{eqn_resolvent_fix_point_iteration},
\end{equation}
we can find a fixed point of $(\bbD+\ccalT)^{-1}\bbD$ and consequently a root of the operator $\ccalT$. Later in Section \ref{sec:results}, we show that the operator $(\bbD+\ccalT)^{-1}$ is well-defined (Lemma \ref{lemma:operators_prop}) and prove by following the update in \eqref{eqn_resolvent_fix_point_iteration}, the iterates converge to a fixed point of  $(\bbD+\ccalT)^{-1}\bbD$. We would like to highlight that the update in \eqref{eqn_resolvent_fix_point_iteration} can be interpreted as performing a proximal-point update, and for this reason we refer to our proposed method as Decentralized Proximal Point for Saddle Point problems (DPPSP).

To ensure that the update in \eqref{eqn_resolvent_fix_point_iteration} can be implemented in a distributed way we define $\bbD$ as
\begin{equation}
\bbD \triangleq  \begin{bmatrix}
\IB & \UB \\
\UB & \IB
\end{bmatrix}.
\label{eq:def_D}
\end{equation}
To check that $\bbD$ is positive definite note that, based on Schur complement, $\bbD\succ \bb0$  holds if and only if $\bbI-\bbU^2= \hbW\succ \bb0$, which is satisfied based on the last condition in \eqref{w_conditions}. 
The operator $(\bbD+\ccalT)^{-1}$ can be implemented in a distributed fashion, as it can be simplified to
\begin{equation}
(\bbD+\ccalT)^{-1} = \begin{bmatrix}
(\alpha[\BM + \RM]+\bbI)^{-1} & (\alpha[\BM + \RM]+\bbI)^{-1}\UB \\
\bb0 & \bbI
\end{bmatrix}.
\label{eq:def_D}
\end{equation}
The operator $(\alpha[\BM + \RM]+\bbI)^{-1}$ has a block diagonal structure and can be computed locally. Moreover, the operator $\bbU$ is graph sparse and can be implemented by exchanging information among neighboring nodes. Therefore, DPPSP is a fully decentralized method.
Now we proceed to formally state the local updates of nodes for implementing DPPSP. By premultiplying both sides of $\eqref{eqn_resolvent_fix_point_iteration}$ by $(\bbD+\ccalT)$ and plugging in the definitions of $\bbv$, $\bbD$, and $ \TM$ we obtain 
\begin{align}
\bbz^{t+1} + \alpha[\BM(\bbz^{t+1}) + \RM(\bbz^{t+1})] &= \bbz^t  + \bbU\bbq^t - 2\bbU\bbq^{t+1}, \label{eqn_derivation_I}\\
\bbq^{t+1} &= \bbU\bbz^t + \bbq^t. \label{eqn_derivation_II}
\end{align}
Substitute $\bbq^{t+1}$ in \eqref{eqn_derivation_I} by its in \eqref{eqn_derivation_II} and use the definition $\bbU^2=\bbI-\hbW$ to obtain
\begin{align}
\bbz^{t+1} + \alpha[\BM(\bbz^{t+1}) + \RM(\bbz^{t+1})] &=(2\hbW-\bbI) \bbz^t  - \bbU\bbq^t, \label{eqn_derivation_a}\\
\bbq^{t+1} &= \bbU\bbz^t + \bbq^t. \label{eqn_derivation_b}
\end{align}
By subtracting two consecutive updates of $\bbz$ we can eliminate $\bbq$ from the update of $\bbz$ and obtain 
\begin{align}\label{final_update}
\bbz^{t+1} & + \alpha[\BM(\bbz^{t+1})  + \RM(\bbz^{t+1})]  \nonumber \\
& =2\hbW \bbz^t -\hbW \bbz^{t-1}  + \alpha[\BM(\bbz^{t}) + \RM(\bbz^{t})],
\end{align}
for $t\geq 1$. By setting $\bbq^0 = \zeroB$, the update for $t = 0$ is given by 
$ \bbz^1 + \alpha[{\BM}(\bbz^{1}) + \RM(\bbz^{1})] =(2\hbW - \IB)\bbz^{0}$. Note as the mixing matrix $\hbW $ has the sparsity pattern of the graph, computation of $2\hbW \bbz^t$ and $\hbW \bbz^{t-1}$ in \eqref{final_update} can be done in a distributed manner. Further, if the operator $(\ccalI+\alpha(\BM + \RM))^{-1} $ is well-defined, then it can be implemented in a distributed fashion as both $\BM$ and $\RM$ have a block diagonal structure. To better illustrate these observations, we would like to highlight that the local version of the update in \eqref{final_update} at node $n$ is given by
\begin{align}\label{final_update_local}
&\bbz_n^{t+1} + \alpha[\BM_n(\bbz_n^{t+1}) + \RM_n(\bbz_n^{t+1})]  \nonumber \\
&=  \sum_{m\in \ccalN_n}w_{nm}(2\bbz_m^t - \bbz_m^{t-1}  )+ \alpha[\BM(\bbz_n^{t}) + \RM(\bbz_n^{t})].
\end{align}
For the first iterate, the update is $\bbz_n^{1} + \alpha[\BM_n(\bbz_n^{1}) + \RM_n(\bbz_n^{1})] 
=  \sum_{m\in \ccalN_n}(2w_{nm}-1)\bbz_m^0$. Hence, node $n$ can perform its update in \eqref{final_update_local}, by computing the proximal operator for $\ccalB_n+\ccalR_n$ or in other words by computing the local operator $(\ccalI+\alpha(\ccalB_n+\ccalR_n))^{-1}$. The steps of the DPPSP method are summarized in Algorithm \ref{alg_main}. 

\begin{algorithm}[t]
	\caption{DPPSP at node $n$}
	\begin{algorithmic}[1]
		\REQUIRE initial iterate $\bbz_n^0$, step size $\alpha$, weights $w_{nm}$ for $m \in \NM_n$;
		\FOR{$t = 0,\dots,T-1$}
		\STATE Exchange variable $\zB_n^t$ with neighboring nodes $m \in \NM_n$;
		\IF{$t=0$}
		\STATE $ \bbz^{t+1} =(\bbI+ \alpha(\BM_n + \RM_n) )^{-1}\big(\sum_{m\in \ccalN_n}(2w_{nm}-1)\bbz_m^t \big)$;
		\ELSE
		\STATE $\bbz_n^{t+1} 
=  (\bbI+ \alpha(\BM_n + \RM_n) )^{-1}\big(\sum_{m\in \ccalN_n}w_{nm}(2\bbz_m^t - \bbz_m^{t-1}  )+ \alpha[\BM(\bbz_n^{t}) + \RM(\bbz_n^{t})]\big)$;
		\ENDIF
		\ENDFOR
	\end{algorithmic}
	\label{alg_main}
\end{algorithm}

%% file: 5-theory.tex

\section{Theoretical Results}\label{sec:results}

In this section, we first show that the inverse operators $(\bbD+\ccalT)^{-1}$ and $(\bbI +\alpha (\ccalB+\ccalR))^{-1}$ are well-defined for a properly chosen $\alpha$. Then, we characterize convergence properties of DPPSP.

\begin{lemma}\label{lemma:operators_prop}
Consider the definitions of the operator $\ccalT$ in \eqref{eqn_augmented_operator} and the matrix $\bbD$ in~\eqref{eq:def_D}. Then,
\vspace{-2mm}
\begin{enumerate}
\item[(i)] The operator $\ccalB+\ccalR$ is $\rho$-weakly monotone.
\vspace{-2mm}
\item[(ii)]  If we set $\alpha<\rho^{-1}$, then the operator $(\bbI +\alpha (\ccalB+\ccalR))^{-1}$ is well-defined.
\vspace{-2mm}
\item[(iii)]  If we choose the stepsize $\alpha $ such that  $\alpha\leq ({1-(1-\lambda_{min}(\bbW))^{1/2}})/({2\rho})$, then the operator $\bbD+\ccalT$ is $(\lambda_{min}(\bbW)/4)$-strongly monotone.
\end{enumerate}
\end{lemma}


The results in Lemma~\ref{lemma:operators_prop} show that if the stepsize $\alpha$ is properly chosen then both operators $(\bbD+\ccalT)^{-1}$ and $(\bbI +\alpha (\ccalB+\ccalR))^{-1}$ are well-defined. Therefore, the update rules in \eqref{eqn_resolvent_fix_point_iteration} and \eqref{final_update} for the DPPSP method are well-defined. 

\begin{theorem}\label{thm:approximate_cnvg}
Consider the DPPSP method outlined in Algorithm \ref{alg_main}. Suppose the conditions in Assumption~\ref{ass:weak_con}-\ref{ass:compact} are satisfied and the stepsize is chosen such that $\alpha \leq 1/(2\rho)$. If we run DPPSP for $T$ iterations and choose one of the iterates $s$ uniformly at random form $1,\dots T$ then we have
\begin{align}
	&\mathbb{E}_{s} \left[\left \|\sum_{n=1}^N\ccalB_n(\bbz_n^{s})+\ccalR_n(\bbz_n^{s}) \right\| \right] \nonumber \\
	&\leq \frac{1}{\alpha} \sqrt{\frac{N}{\lambda_{\min}(\hbW)}} \left(\frac{\|\bbphi^0 - \bbphi^*\|_\MB}{\sqrt{T}}  + \sqrt{2\alpha \rho N} D \right),
\end{align}
and 
\begin{align}
	\mathbb{E}_{s} [ \|\bbU\bbz^{s+1}\|]\leq \frac{\|\bbphi^0 - \bbphi^*\|_\MB}{\sqrt{T}}  + \sqrt{2\alpha \rho N} D,
	\end{align}
	where {$\|\bbphi^0 - \bbphi^*\|_\MB= {\| \bbz_0 - \bbz^*\|_{\hbW} + \| \bbU \bbz_0 - \bbq^* \|}$, and $\bbz^*$ and $\bbq^*$ are defined in \eqref{eqn_optimality}. }
\end{theorem}

The result in Theorem \ref{thm:approximate_cnvg} shows that after $T$ iterations if we choose one of the iterates uniformly at random, then in expectation the iterates satisfy first-order optimality condition up to an error of $\mathcal{O}(\sqrt{\rho} +{1}/{\sqrt{T}})$, i.e., $\mathbb{E}_{s} [\|\sum_{n=1}^N\ccalB_n(\bbz_n^{s})+\ccalR_n(\bbz_n^{s}) \|]\leq \mathcal{O}(\sqrt{\rho} +{1}/{\sqrt{T}})$, and the expected consensus error is also of $\mathcal{O}(\sqrt{\rho} +{1}/{\sqrt{T}})$, i.e., $\mathbb{E}_{s} [ \|\bbU\bbz^{s+1}\|]\leq \mathcal{O}(\sqrt{\rho} +{1}/{\sqrt{T}})$. These results show that the first-order optimality condition gap and the consensus constraint violation sublinearly at a rate of $\mathcal{O}({1}/{\sqrt{T}})$ converges to a neighborhood that has a radius of $\mathcal{O}(\sqrt{\rho})$. Hence, as the function $f$ becomes closer to a convex-concave function, i.e., $\rho$ becomes smaller, the accuracy of our results become better. In particular, when $f$ is convex-concave, i.e., $\rho=0$, the iterates generated by DPPSP find a point that as $\eps$-optimality gap and $\eps$-consensus error after at most $\mathcal{O}(1/\eps^2)$ iterations.
To achieve any arbitrary first-order stationary point we need to add the following assumption. 

\begin{assumption}\label{MVI_assumption}
(Minty Variational Inequality [MVI]): There exists a point $\bbz^*$ such that for every local operator $\ccalB_n$,  $\ccalB_n(\bbz)^\top (\bbz-\bbz^*)\geq 0$ for any $\bbz\in \ccalZ$.
\end{assumption}
This is a standard assumption which is made for non-convex optimization problems (see \cite{landang}, \cite{lin2018solving} for more details). 
This assumption is clearly satisfied when the operator is monotone (which corresponds to the function being convex-concave i.e. $\rho = 0$). Pseudo-monotone operators i.e. operators which satisfy the property: 
\begin{align}
\langle F(y),x - y \rangle \geq 0 \implies \langle F(x),x-y \rangle \geq 0,
\end{align}
also satisfy this condition. There are also other classes of operators which satisfy the Minty's VI condition. For example, consider the following operator which is neither monotone nor pseudo-monotone (example from \cite{landang})
\begin{align}
F(x) = 
\begin{cases} 
      0 & x\ = x_0  \\
      \geq 0 & x \geq x_0 \\
      \leq 0 & x \leq x_0 
   \end{cases}
\end{align}
The solution to this problem is $x^* = x_0$, and at this point it can be easily seen that the Minty VI condition is satisfied.

Once we have this assumption, we can show exact convergence of DPPSP to a solution of Problem \eqref{eqn_decentralized_consensus_root_monotone_operator} as shown in the following theorem.

%
%

	\begin{figure*}[t!]
		\begin{subfigure}[b]{0.5\textwidth}
			\centering
			\includegraphics[width=\textwidth]{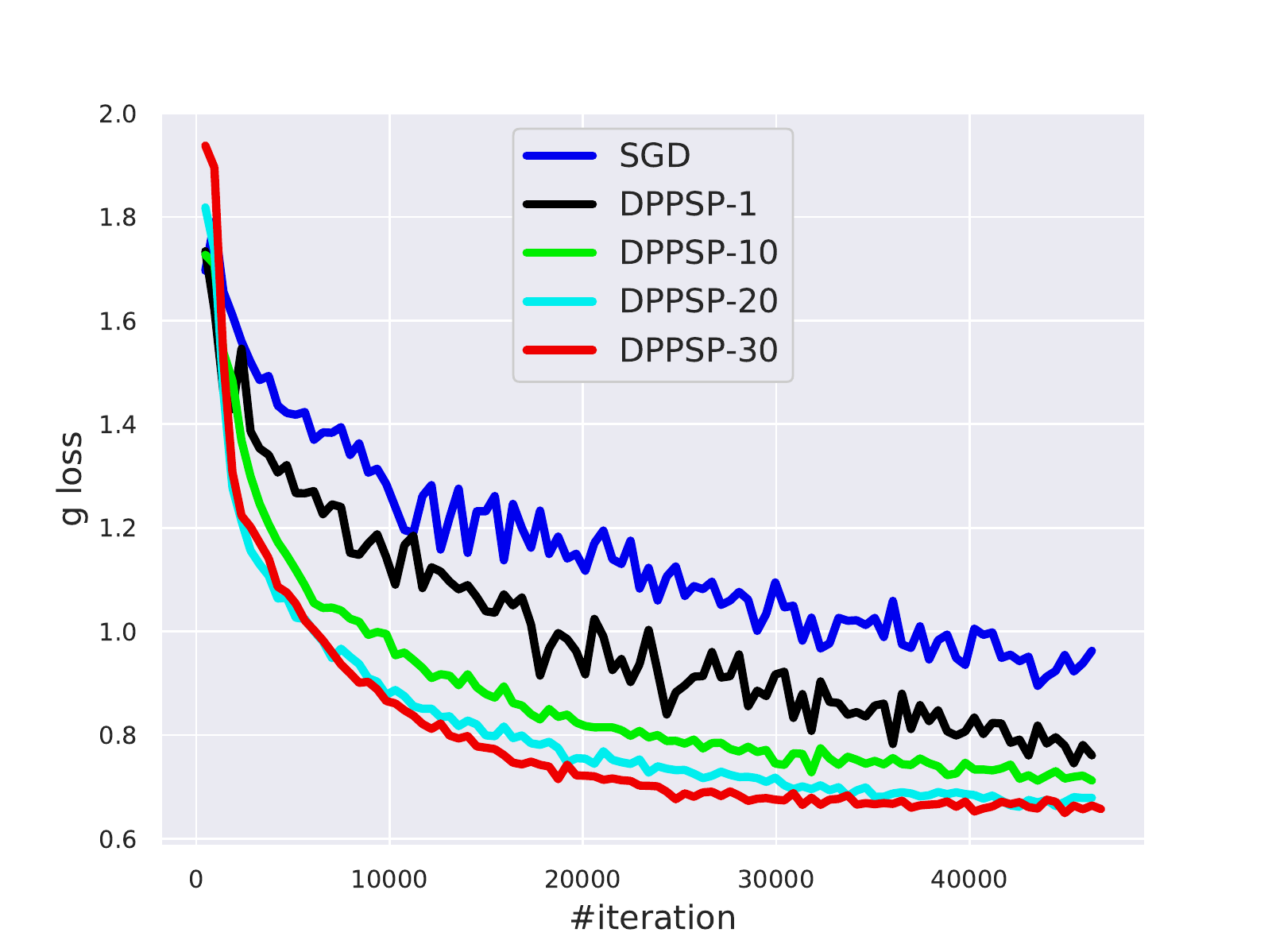}
			\caption{Generator Loss}
			\label{ex_sgd}
		\end{subfigure}
		\hfill
		\begin{subfigure}[b]{0.5\textwidth}
			\centering
			\includegraphics[width=\textwidth]{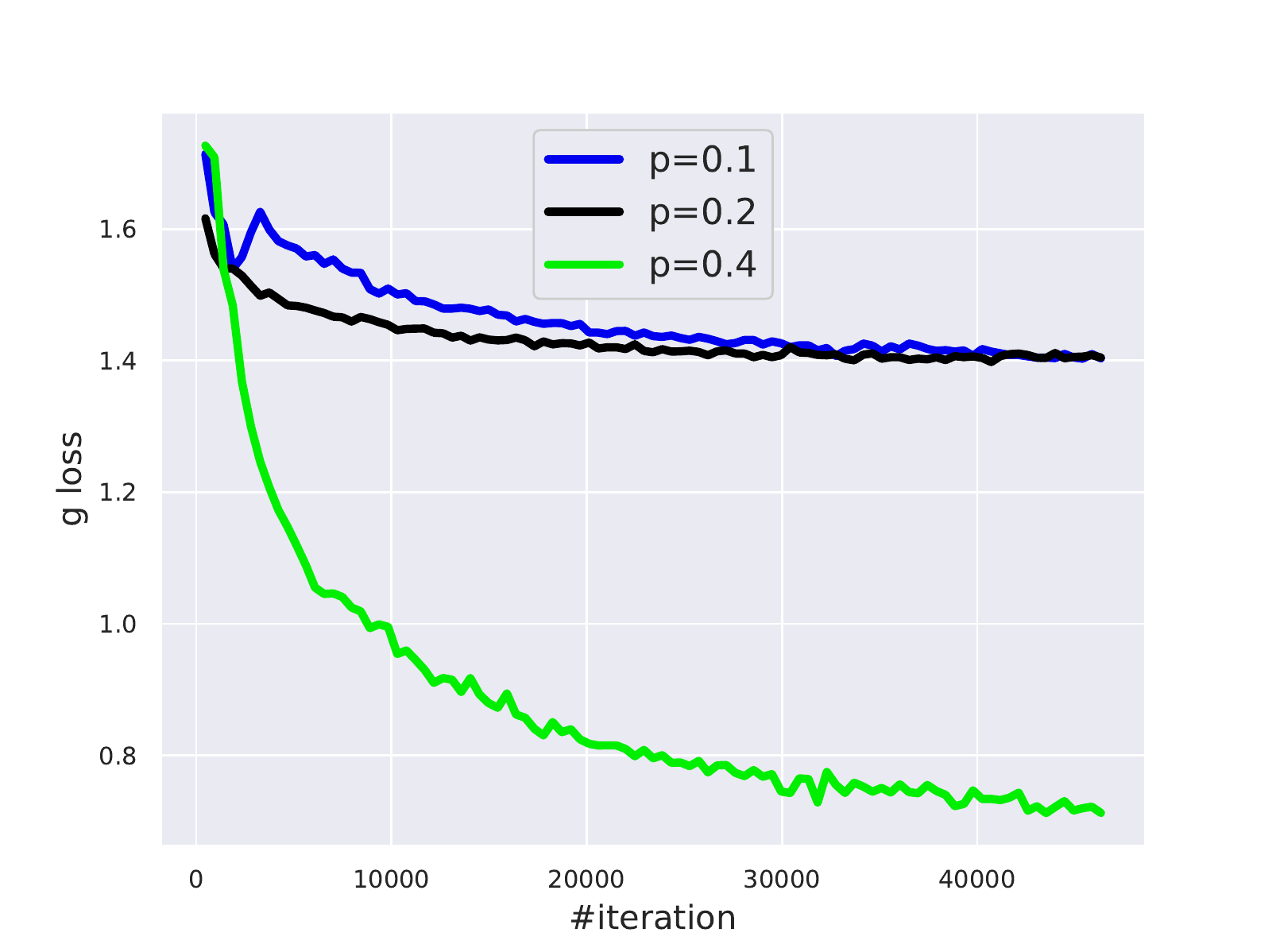}
			\caption{Different Sparsity}
			\label{ex_p}
		\end{subfigure}
		\hfill
	\caption{Performances of DPPSP with varying number of nodes and different graph connectivity on the MNIST dataset. Single-node DPPSP method outperforms SGD. The advantage of DPPSP is more significant when we have more processing nodes. Note that the number of samples used per iteration are identical for all the considered settings and hence the the number of iterations is proportional to the wall-time of each method. }
	\end{figure*}
\begin{theorem}
\label{thm:convergence_DPPSP}
Consider the DPPSP method outlined in Algorithm \ref{alg_main}. Suppose the conditions in Assumption~\ref{ass:weak_con}-\ref{MVI_assumption} are satisfied and the stepsize is chosen such that $\alpha = 1/(2\rho)$. If we run DPPSP for $T$ iterations and choose one of the iterates $s$ uniformly at random form $1,\dots, T$ then we have
\begin{align}
	\mathbb{E}_{s} \left[\left \|\sum_{n=1}^N\ccalB_n(\bbz_n^{s+1})+\ccalR_n(\bbz_n^{s+1}) \right\|\right]
	&\leq \frac{ND}{\alpha \sqrt{T}},\nonumber \\
	\mathbb{E}_{s} \left[ \|\bbU\bbz^{s}\|\right]
	&\leq \frac{\sqrt{N}D}{\sqrt{T}}.
	\end{align}
\end{theorem}

The result in Theorem~\ref{thm:convergence_DPPSP} shows that once the MVI assumption holds, the iterates generated by DPPSP can achieve any arbitrary $\eps$ accuracy. In particular, they find a solution with an $\eps$-first-order optimality gap and an $\eps$-consensus error after at most $\mathcal{O}(1/\eps^2)$ iterations.

%% file: 7-experiments.tex

\section{Numerical Experiments}
\label{sec:exp}

\begin{table*}[t]
	\centering
	\caption{Fake images produced by generators that are trained using SGD and variants of DPPSP.}\label{tbl:samples}
	\begin{tabular}{ | c | c  c  c  c  c| }
		\hline
		iterations& 2340 & 4680 & 9360 & 18720 & 28080\\ 
		\hline
		SGD
		&
		\begin{minipage}{.14\textwidth}
			\includegraphics[width=\linewidth, height=18mm]{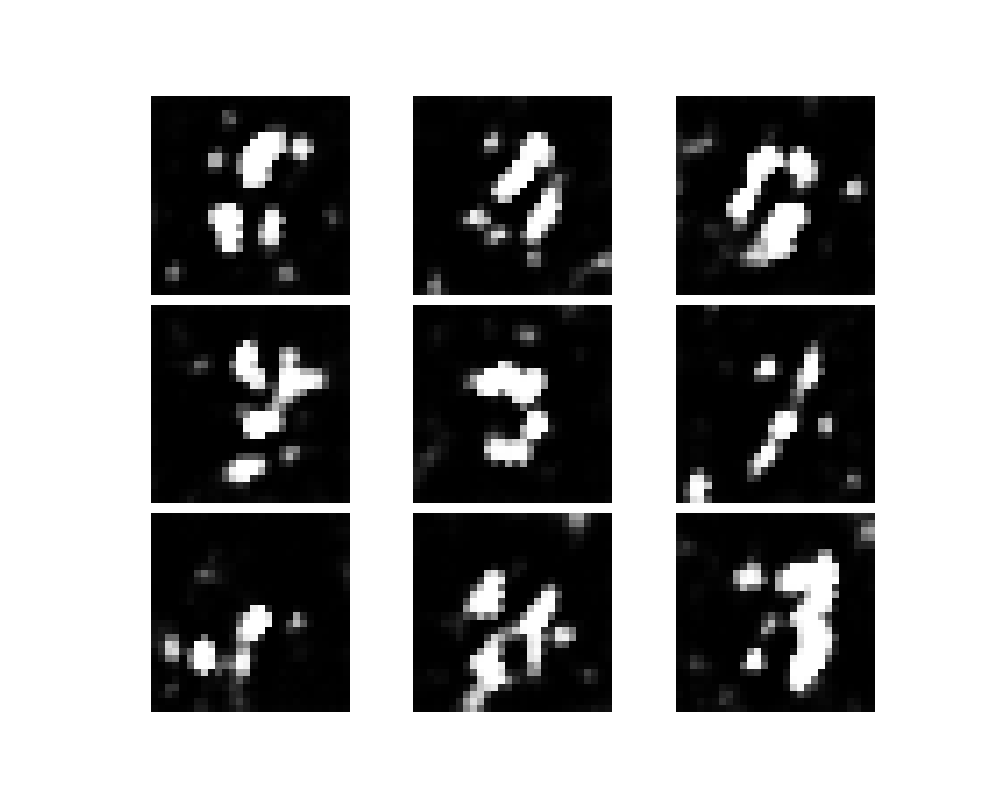}
		\end{minipage}
		&
		\begin{minipage}{.14\textwidth}
			\includegraphics[width=\linewidth, height=18mm]{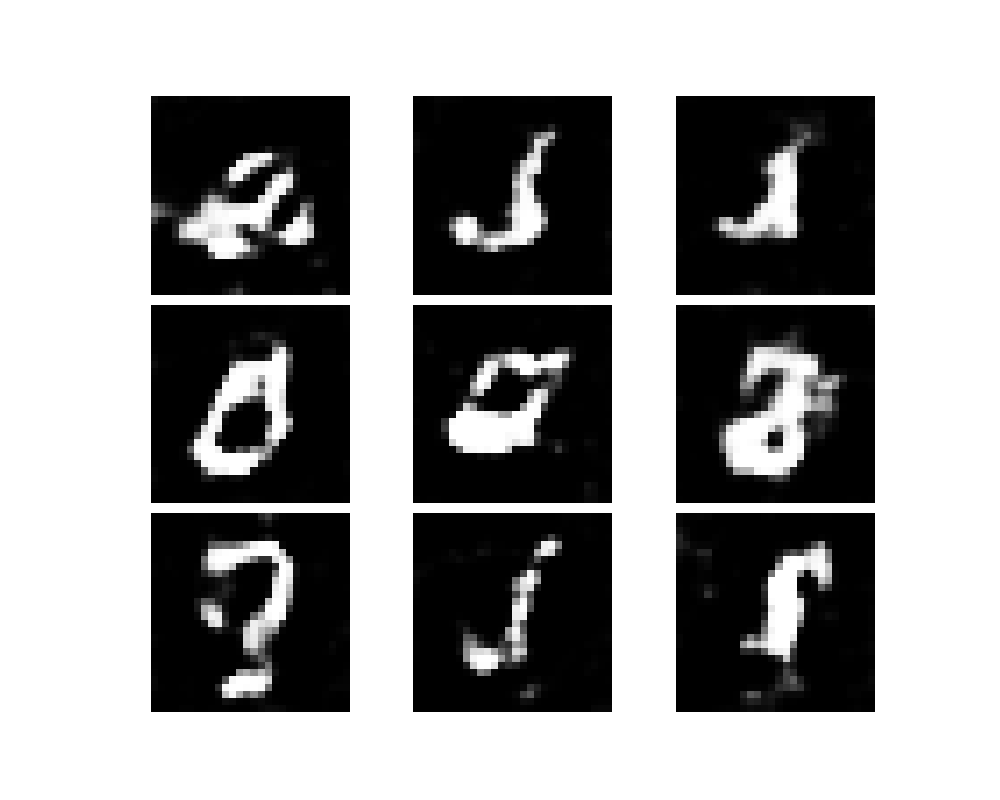}
		\end{minipage}
		&
		\begin{minipage}{.14\textwidth}
			\includegraphics[width=\linewidth, height=18mm]{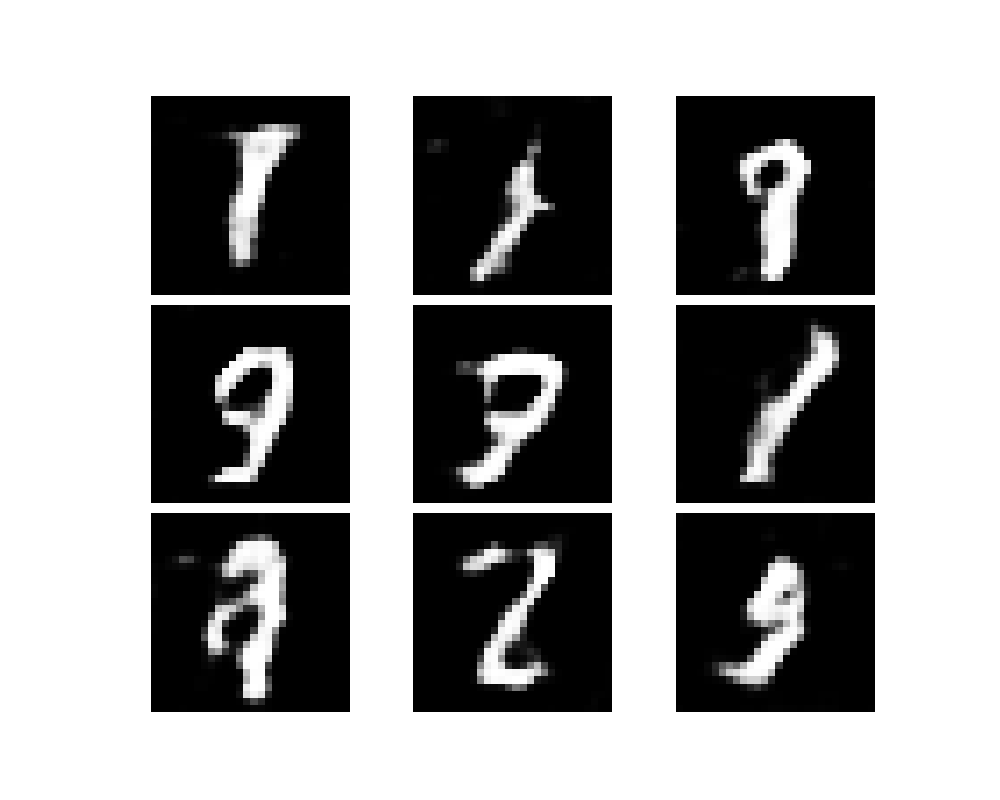}
		\end{minipage}
		&
		\begin{minipage}{.14\textwidth}
			\includegraphics[width=\linewidth, height=18mm]{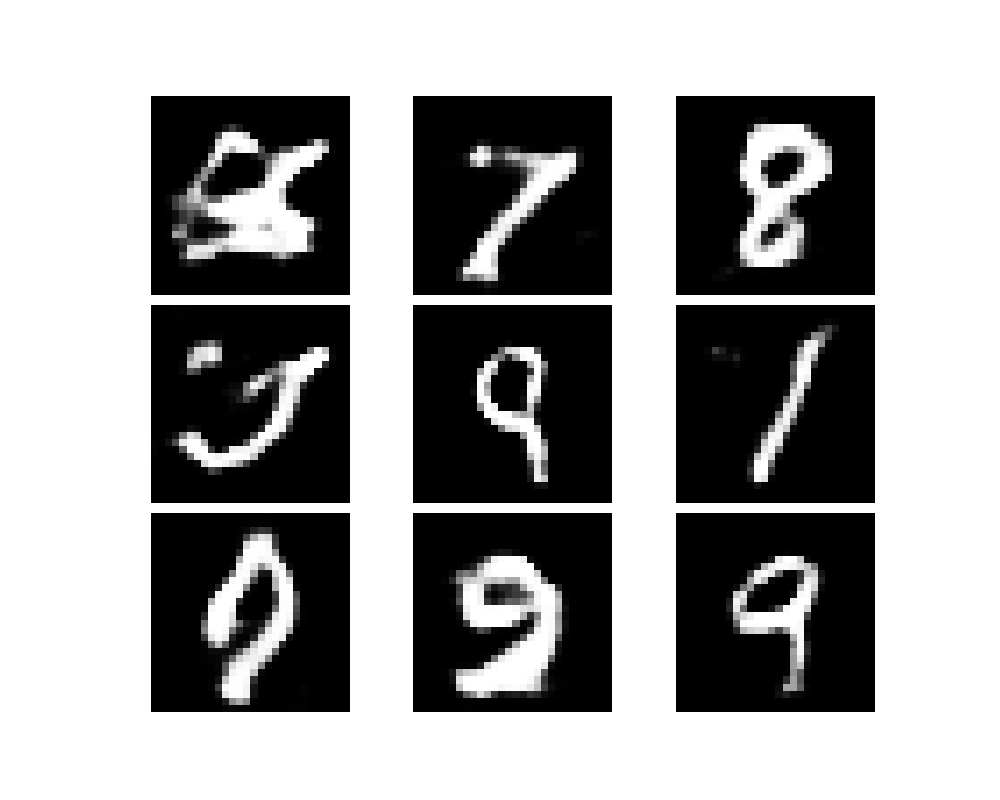}
		\end{minipage}
		&
		\begin{minipage}{.14\textwidth}
			\includegraphics[width=\linewidth, height=18mm]{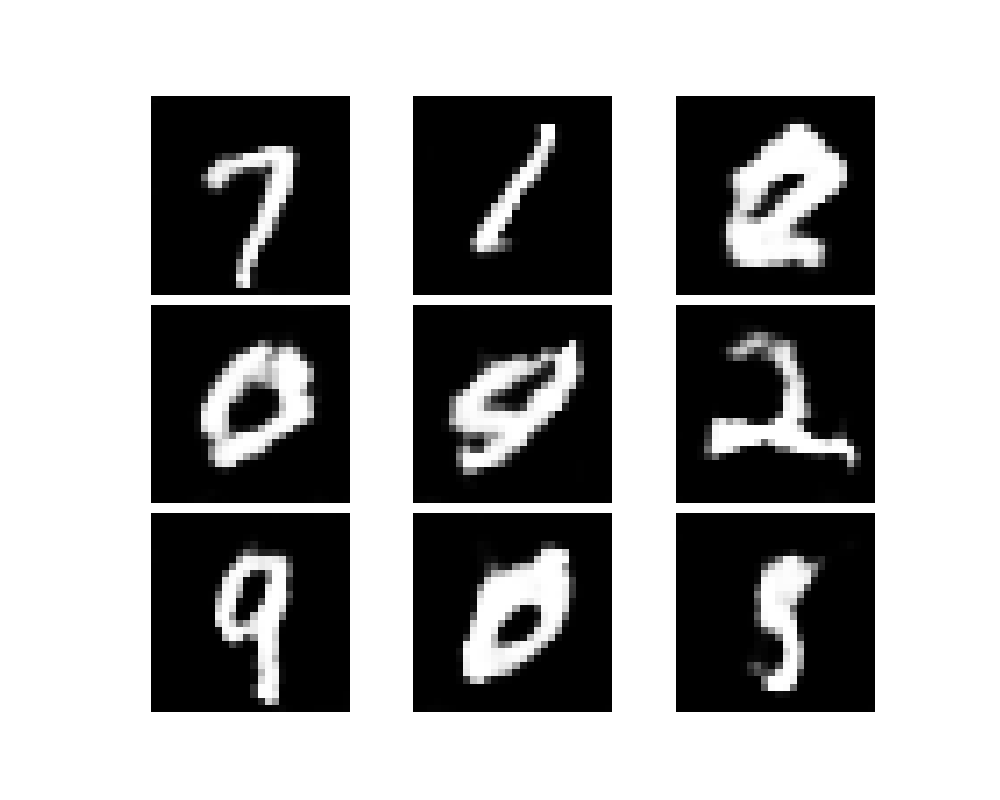}
		\end{minipage}
		\\ 
		DPPSP-1
		&
		\begin{minipage}{.14\textwidth}
			\includegraphics[width=\linewidth, height=18mm]{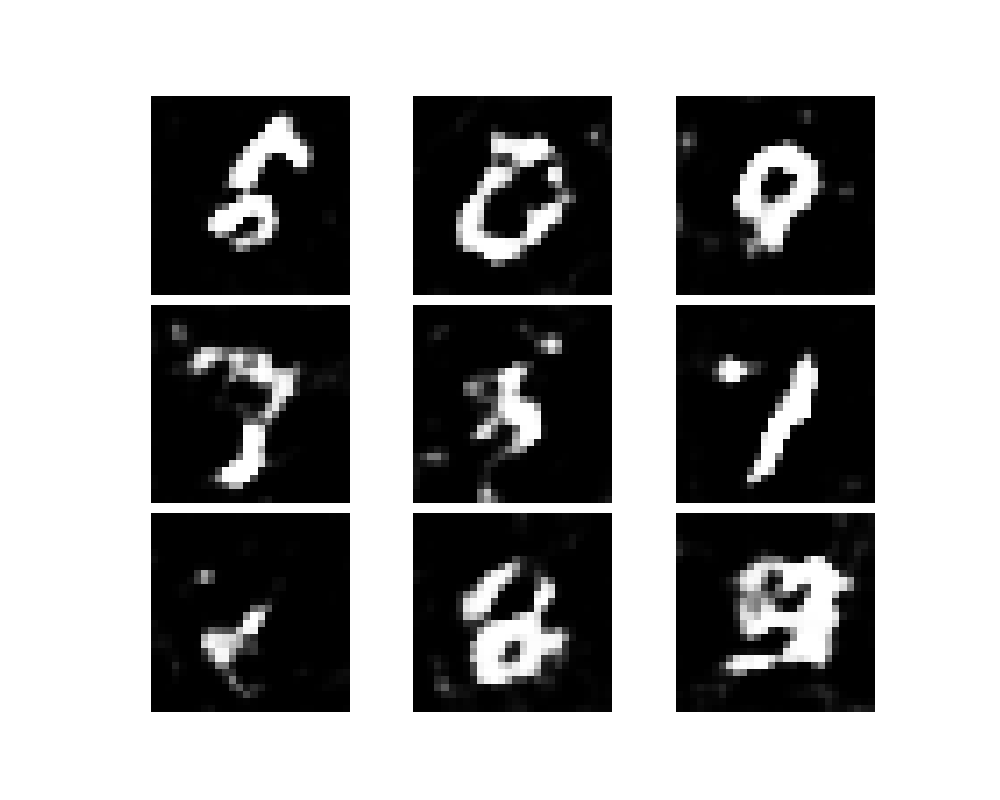}
		\end{minipage}
		&
		\begin{minipage}{.14\textwidth}
			\includegraphics[width=\linewidth, height=18mm]{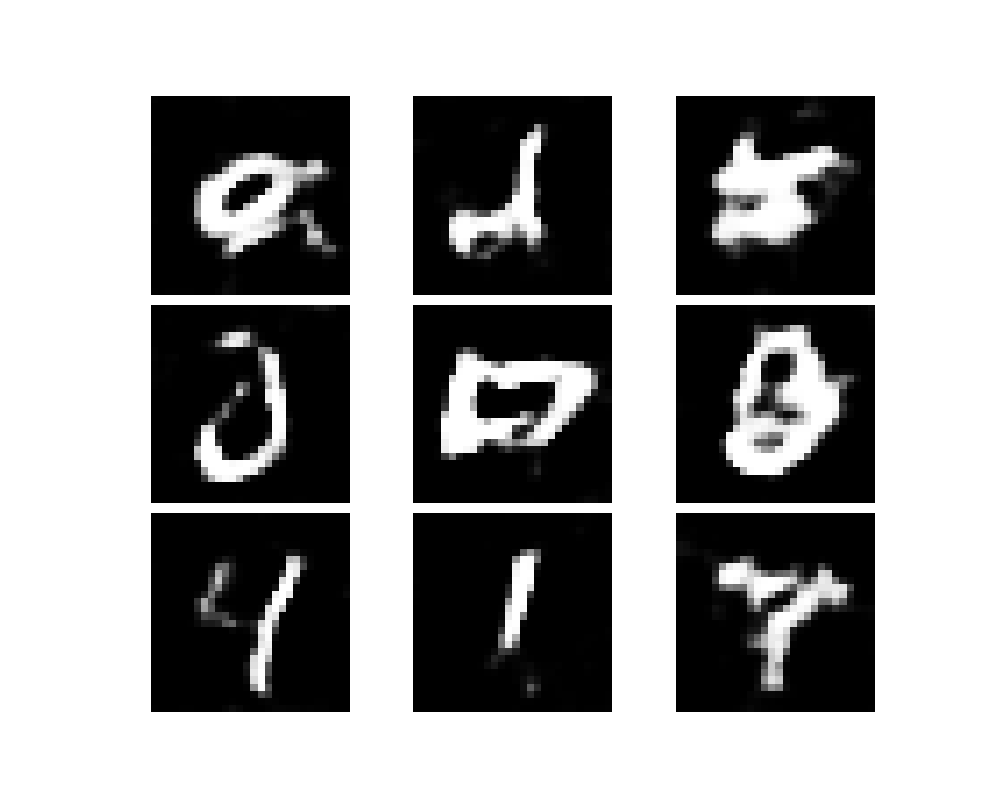}
		\end{minipage}
		&
		\begin{minipage}{.14\textwidth}
			\includegraphics[width=\linewidth, height=18mm]{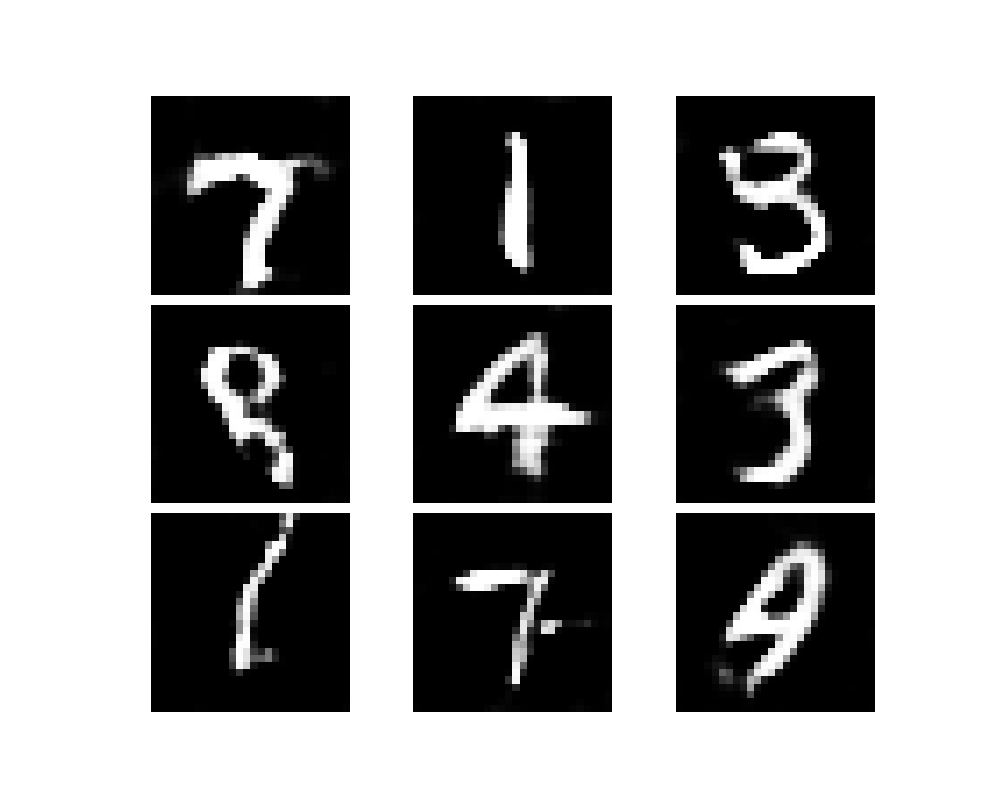}
		\end{minipage}
		&
		\begin{minipage}{.14\textwidth}
			\includegraphics[width=\linewidth, height=18mm]{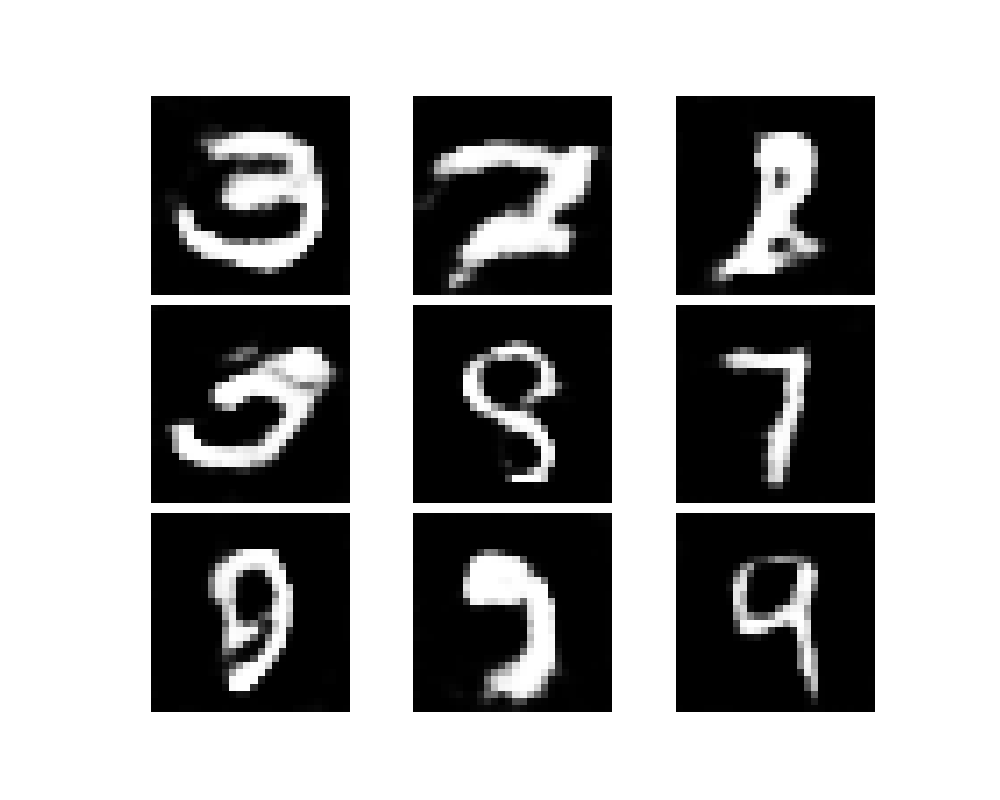}
		\end{minipage}
		&
		\begin{minipage}{.14\textwidth}
			\includegraphics[width=\linewidth, height=18mm]{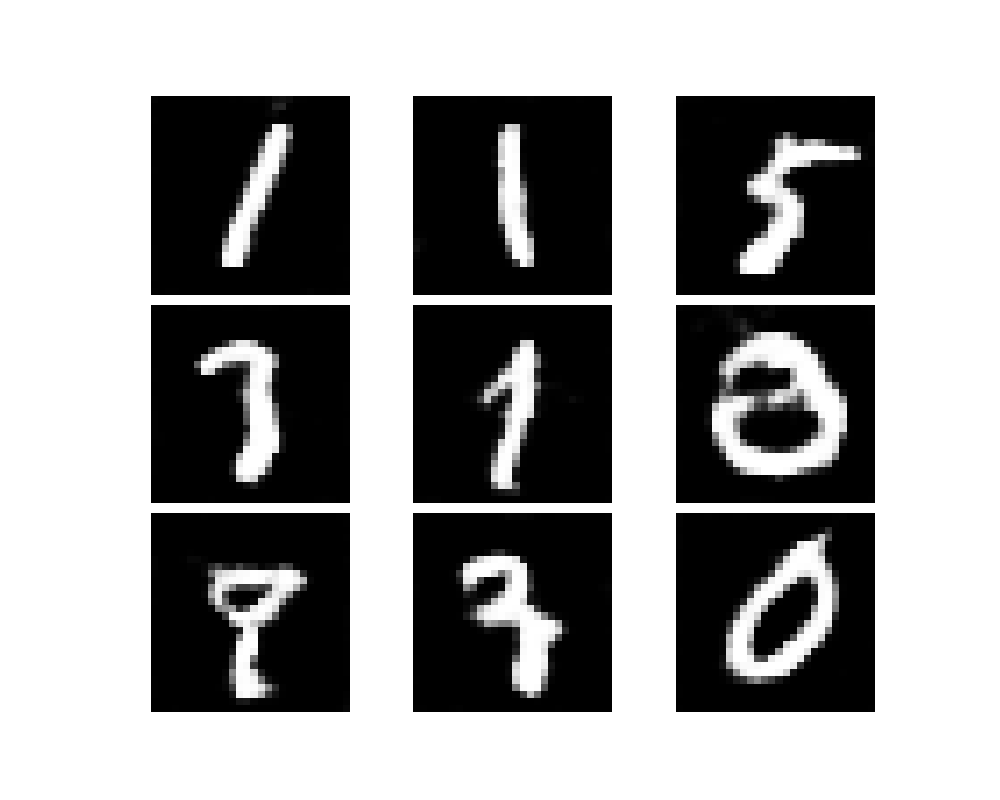}
		\end{minipage}
		\\ 
		DPPSP-10
		&
		\begin{minipage}{.14\textwidth}
			\includegraphics[width=\linewidth, height=18mm]{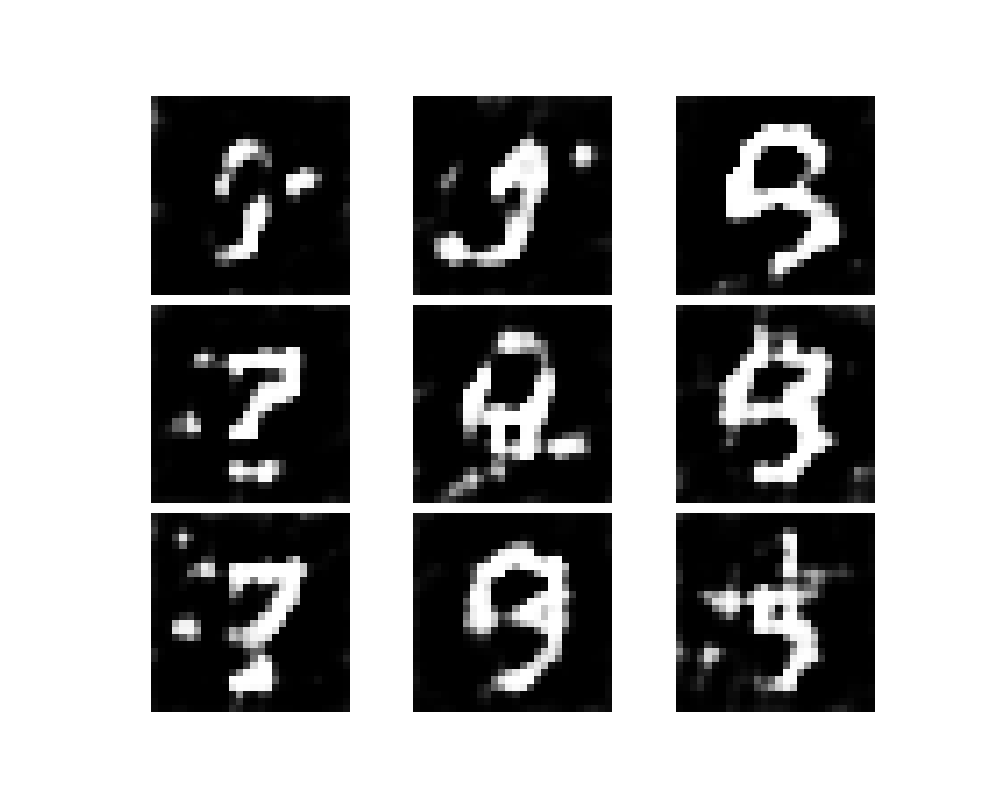}
		\end{minipage}
		&
		\begin{minipage}{.14\textwidth}
			\includegraphics[width=\linewidth, height=18mm]{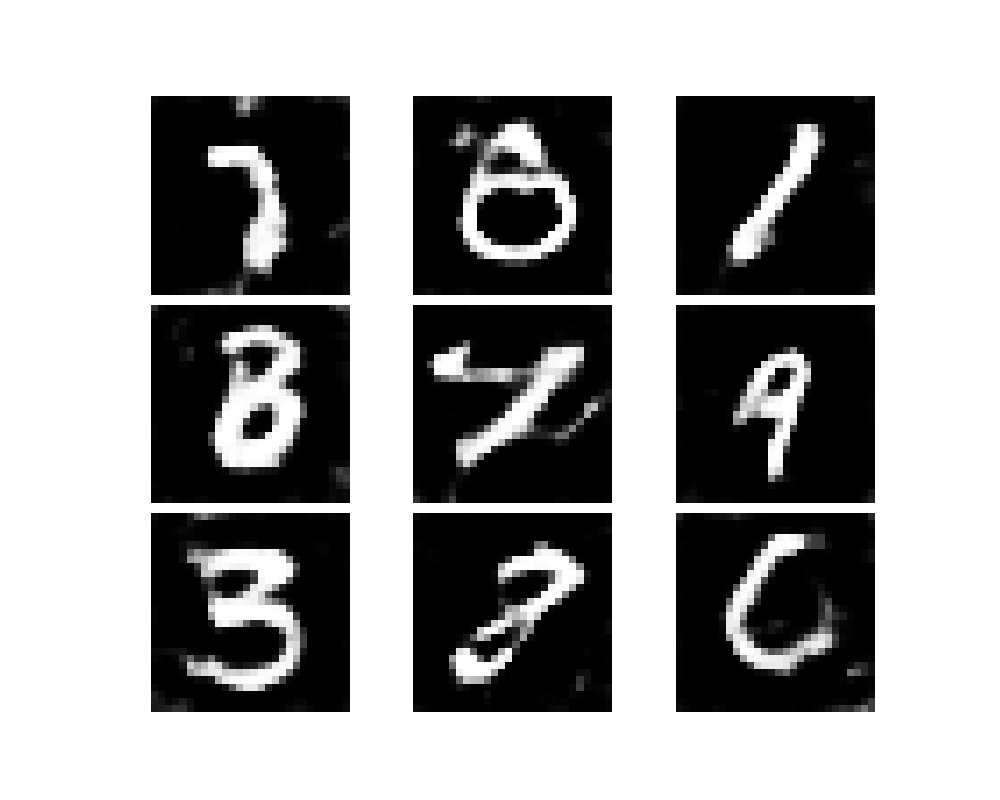}
		\end{minipage}
		&
		\begin{minipage}{.14\textwidth}
			\includegraphics[width=\linewidth, height=18mm]{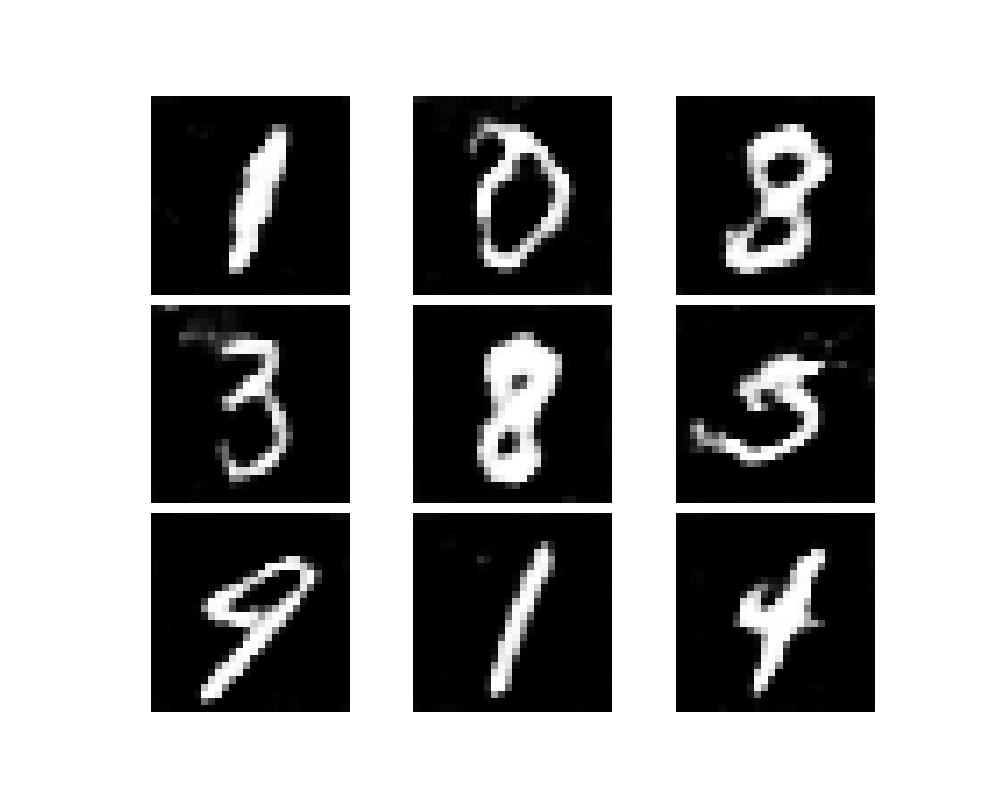}
		\end{minipage}
		&
		\begin{minipage}{.14\textwidth}
			\includegraphics[width=\linewidth, height=18mm]{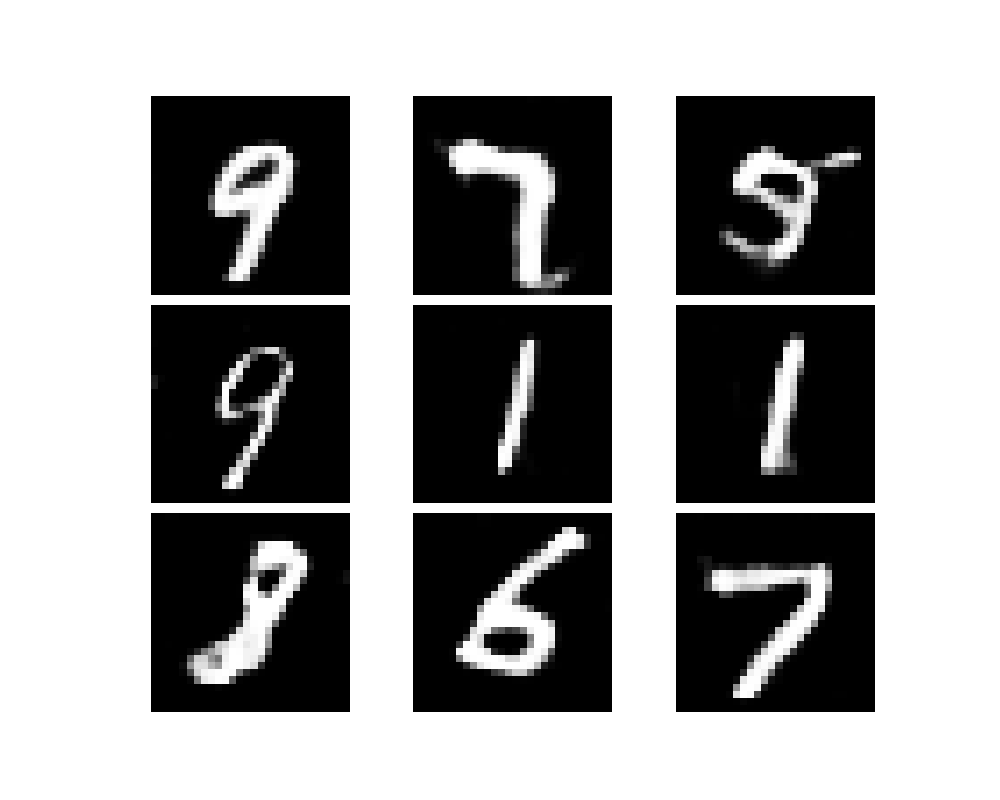}
		\end{minipage}
		&
		\begin{minipage}{.14\textwidth}
			\includegraphics[width=\linewidth, height=18mm]{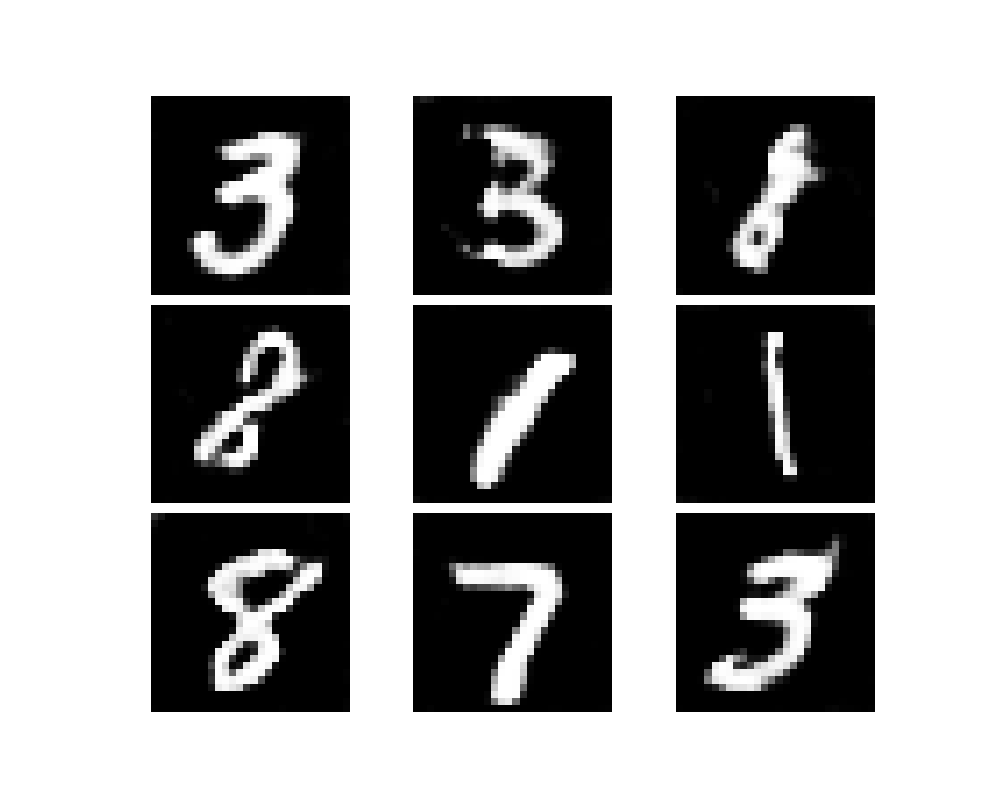}
		\end{minipage}
		\\ 
		DPPSP-20
		&
		\begin{minipage}{.103\textwidth}
			\includegraphics[width=\linewidth, height=15mm]{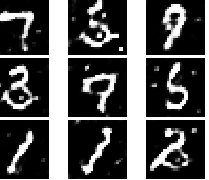}
		\end{minipage}
		&
		\begin{minipage}{.103\textwidth}
			\includegraphics[width=\linewidth, height=15mm]{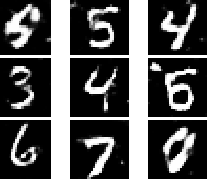}
		\end{minipage}
		&
		\begin{minipage}{.103\textwidth}
			\includegraphics[width=\linewidth, height=15mm]{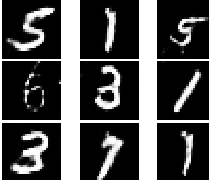}
		\end{minipage}
		&
		\begin{minipage}{.103\textwidth}
			\includegraphics[width=\linewidth, height=15mm]{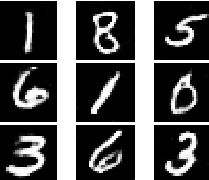}
		\end{minipage}
		&
		\begin{minipage}{.103\textwidth}
			\includegraphics[width=\linewidth, height=15mm]{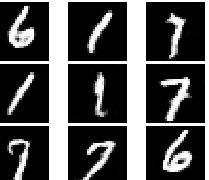}
		\end{minipage}
		\\[23pt]
		DPPSP-30
		&
		\begin{minipage}{.103\textwidth}
			\includegraphics[width=\linewidth, height=15mm]{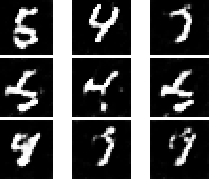}
		\end{minipage}
		&
		\begin{minipage}{.103\textwidth}
			\includegraphics[width=\linewidth, height=15mm]{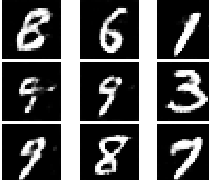}
		\end{minipage}
		&
		\begin{minipage}{.103\textwidth}
			\includegraphics[width=\linewidth, height=15mm]{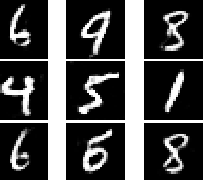}
		\end{minipage}
		&
		\begin{minipage}{.103\textwidth}
			\includegraphics[width=\linewidth, height=15mm]{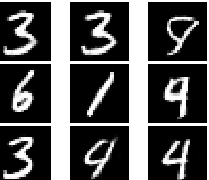}
		\end{minipage}
		&
		\begin{minipage}{.103\textwidth}
			\includegraphics[width=\linewidth, height=15mm]{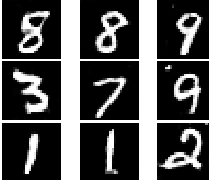}
		\end{minipage}
		\\[23pt] \hline
	\end{tabular}
\end{table*}
In this section, we show the empirical performance of DPPSP in the context of GAN training\footnote{Our code can be downloaded from https://1drv.ms/u/s!AnXFkPphDb-Ng2cR6ooO6sMY-N8d}. 
We use the number of iterations to measure the computational cost of each algorithm.
The number of samples used per-iteration are identical for all settings and hence the the number of iterations is proportional to the wall-time. 
To assess the capability of generators in training, we run the code downloaded from \cite{cs231n} to produce a pair of ``optimal'' discriminator and generator, which are able to generate high-quality and diversified images.
We use the optimal discriminator to evaluate the performance of the local generators.
The global $g_{loss}$ is defined to be the average of local losses.
%
Note that we do not use the optimal discriminator and optimal generator to guide network training. In all experiments, we assign $L=25$. Our experiments are conducted on MNIST dataset as well as celebA dataset. 

We use SGD as baseline to validate the performance of the proposed DPPSP method.
Concretely, Figure \ref{ex_sgd} shows the (i) the advantage of DPPSP holds over SGD (even when $N=1$) and (ii) how its performance improves when we have more nodes for MNIST dataset.
SGD is run on a single machine which stores all of the training data. Single-node DPPSP has $N=1$ and $\mathbf{W}=1$. It also has access to all data like SGD. For $N>1$, the training data is randomly and equally split to the $N$ machines. We generate graph edges with probability 0.4 and set $\mathbf{W}=\mathbf{I}-{\mathbf{L}}/{\tau}$, where $\mathbf{L}$ is the Laplacian matrix and  $\tau\ge{\lambda_{max}(\mathbf{L})}/{2}$ is a scaling parameter. The graph structure is constant throughout training. We conduct extensive search over hyperparameter $\alpha$ and report best results. It can be seen that the single-node DPPSP method outperforms SGD. The advantage of DPPSP gets clearer when we have more nodes, which shows the scalability of DPPSP. We can also observe that DPPSP becomes more stable when more nodes are available in the network.

\begin{figure}[t!]
	\centering
	\includegraphics[width=0.6\textwidth]{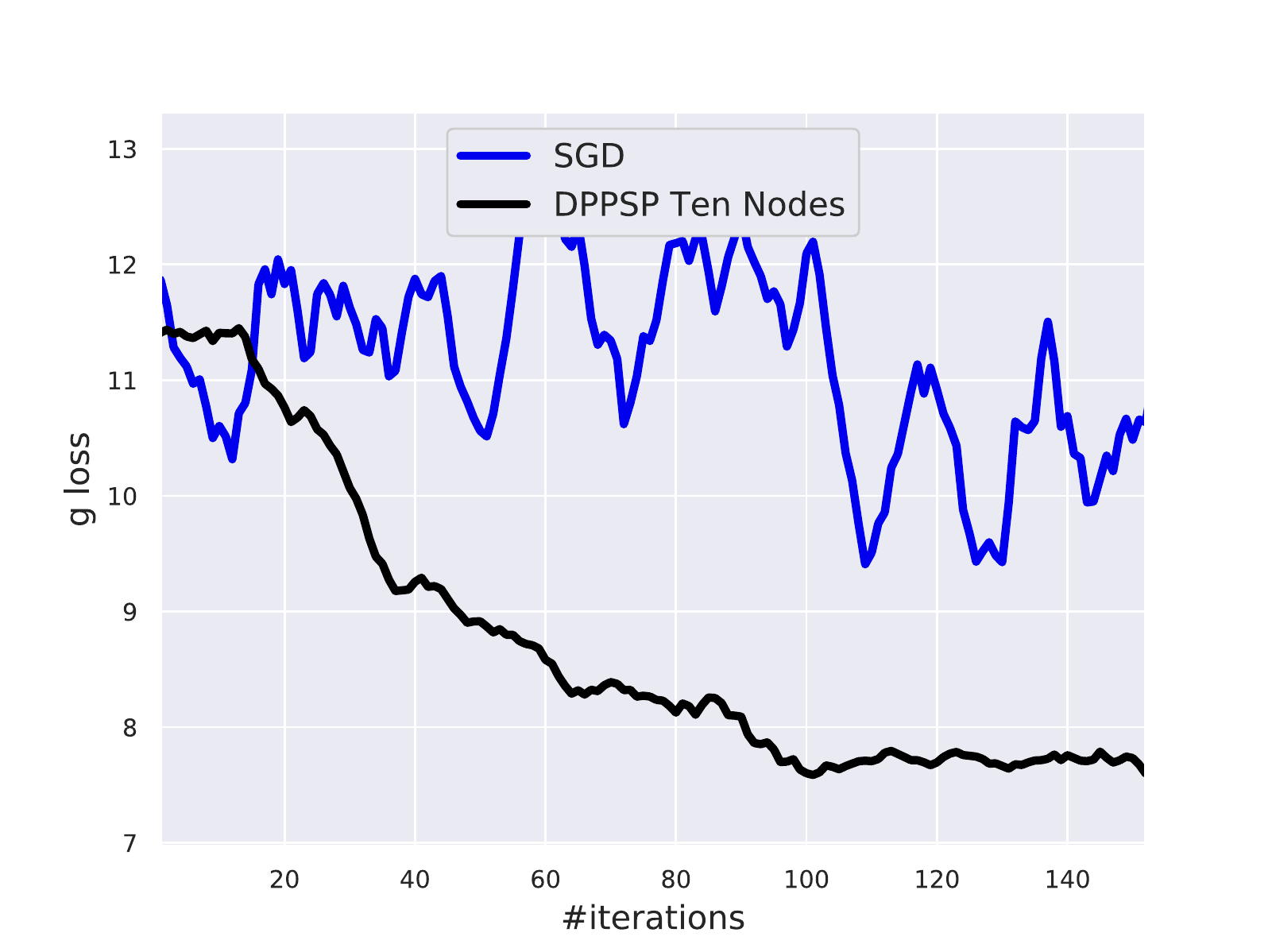}
	\caption{Performance (Generator Loss) of DPPSP with varying number of nodes on the celebA dataset.}
	\label{fig:celebA}
\end{figure}

To show the impact of the graph structure, we vary the edge sparsity and test the proposed DPPSP method on ten-node graphs. The edges are generated randomly with probabilities $0.1$, $0.2$, and $0.4$ with corresponding edge numbers of  $21$, $27$, and $45$, respectively. As we observe in Figure \ref{ex_p}, DPPSP with $p=0.4$ has the best graph connectivity and hence has the best performance.

To better illustrate the advantage of DPPSP over SGD we also compare the images produced by these two methods as time progresses for the MNIST dataset. Table~\ref{tbl:samples} shows images produced during training for SGD and DPPSP with $1$, $10$, $20$, and $30$ nodes. We present the images produced by these methods after $2340$, $4680$, $9360$, $18720$, and $28080$ iterations. As we can observe, after $4680$ iterations, DPPSP-20 and DPPSP-30 already have generated reasonable images, in particular, images of DPPSP-30 have better quality than the ones generated by DPPSP-20. Generators trained by SGD and DPPSP-1 output unsatisfactory samples even after $28080$ iterations. According to these results, increasing the number of processing units not only leads to a smaller loss for the generator but also produces images with higher quality.

We also compare SGD and different variants of DPPSP on the celebA dataset. Figure \ref{fig:celebA} demonstrates the loss of generator corresponding to SGD and DPPSP with 10 processing nodes. Similarly, here we also observe that the generator loss associated with DPPSP is smaller than the one for SGD. Also, the DPPSP algorithm is more stable compared to SGD.

%% file: 6-appendix.tex

\section{Appendix}
\label{sec:supp}

\subsection{Proof of Lemma \ref{lemma:operators_prop}}\label{app:lemma:operators_prop}

First, we show that the operator $\ccalB+\ccalR$ is $\rho$-weakly monotone. To do so, note that 
	\begin{align}
	& \langle   (\ccalB+\ccalR)(\bbz_1)- (\ccalB+\ccalR)(\bbz_2) ,  \bbz_1- \bbz_2\rangle \nonumber\\
	 &\qquad = \langle  \ccalB(\bbz_1)- \ccalB(\bbz_2) ,  \bbz_1- \bbz_2\rangle +\langle   \ccalR(\bbz_1)-\ccalR(\bbz_2) ,  \bbz_1- \bbz_2\rangle \nonumber\\
	 &\qquad \geq  \langle  \ccalB(\bbz_1)- \ccalB(\bbz_2) ,  \bbz_1- \bbz_2\rangle \nonumber\\
	 &\qquad \geq -\rho \| \bbz_1- \bbz_2\|^2
	\end{align}
where the first inequality holds due to the fact that the sets $\ccalX$ and $\ccalY$ are convex and therefore their indicator functions are also convex which implies that the operator $\ccalR$ is monotone, and the second inequality follows from the fact that operator $\ccalB$ is $\rho$-weakly monotno. Therefore, the operator $\ccalB+\ccalR$ is $\rho$-weakly monotone, and if we choose $\alpha<1/\rho$ then $(\ccalI +\alpha (\ccalB+\ccalR))^{-1}$ is well-defined.

Now we proceed to show that the operator $\ccalE=\bbD+\ccalT$ is strongly monotone. Note that we can write the operator $\ccalE$ as
\begin{equation}
\ccalE(\bbv) = \bigg(
\underbrace{\begin{bmatrix}
	\alpha[\BM + \RM]& 0 \\
	0& 0
	\end{bmatrix}}_{\ccalE_1} +
\underbrace{\begin{bmatrix}
	\bbI & 2\UB \\
	\bb0 & \bbI
	\end{bmatrix}}_{\ccalE_2}
\bigg)
\underbrace{\begin{bmatrix}
	\bbz \\
	\bbq
	\end{bmatrix}}_{\bbv},
\label{proof_100}
\end{equation}
Based on the first part of the proof we know that $\alpha (\ccalB+\ccalR)$ is $\alpha \rho$-weakly monotone. Therefore, the operator $\ccalE_1$ is also $\alpha \rho$-weakly monotone. Now we proceed to show that $\ccalE_2$ is strongly monotone. We can write
\begin{align}
&\langle \ccalE_2(\bbv_2) -\ccalE_2(\bbv_2) , \bbv_2 -\bbv_1\rangle\nonumber\\
&\qquad = \langle (\bbz_2+2\bbU\bbq_2)- (\bbz_1+2\bbU\bbq_1), \bbz_2-\bbz_1\rangle
+  \langle \bbq_2- \bbq_1, \bbq_2-\bbq_1\rangle\nonumber\\
&\qquad = \langle \bbz_2-\bbz_1, \bbz_2-\bbz_1\rangle
+2\langle \bbU(\bbq_2-\bbq_1), \bbz_2-\bbz_1\rangle
+  \langle \bbq_2- \bbq_1, \bbq_2-\bbq_1\rangle \nonumber\\
&\qquad \geq  \|\bbz_2-\bbz_1\|^2
-2\|\bbU\| \|\bbq_2-\bbq_1\| \| \bbz_2-\bbz_1\|
+  \|\bbq_2- \bbq_1\|^2 \nonumber\\
& \qquad 
\end{align}
If we define $\beta:=\max\{ |\lambda_{min}( \bbU)|,|\lambda_{max}( \bbU)|\}$, then it can verified that $\|\bbU\|\leq \beta$. Using this inequality we can show that 
\begin{align}\label{proof_800}
\langle \ccalE_2(\bbv_2) -\ccalE_2(\bbv_2) , \bbv_2 -\bbv_1\rangle
&  \geq  \|\bbz_2-\bbz_1\|^2
-2\beta\|\bbq_2-\bbq_1\| \| \bbz_2-\bbz_1\|
+  \|\bbq_2- \bbq_1\|^2 \nonumber\\
&  \geq  \|\bbz_2-\bbz_1\|^2
-\beta\|\bbq_2-\bbq_1\|^2-\beta \| \bbz_2-\bbz_1\|^2
+  \|\bbq_2- \bbq_1\|^2 \nonumber\\
&  =   (1-\beta) \|\bbz_2-\bbz_1\|^2
+  (1-\beta) \|\bbq_2- \bbq_1\|^2 \nonumber\\
& =   (1-\beta) \|\bbv_2-\bbv_1\|^2.
\end{align}
Now note that $\bbU=(\bbI-\hbW)^{1/2}$. Since the eigenvalues of $\bbW$ belong to the interval $(0,1]$ it can be easily verified that $\bbU$ is positive semidefinite. We can further show that $\lambda_{max}( \bbU)= (1-\lambda_{min}(\bbW))^{1/2}$ and $\lambda_{min}( \bbU)= (1-\lambda_{max}(\bbW))^{1/2}$. Therefore, $\beta=(1-\lambda_{min}(\bbW))^{1/2} $. This result 
implies that $\ccalE_2$ is $1-(1-\lambda_{min}(\bbW))^{1/2} $-strongly monotone. Therefore we can write
\begin{align}\label{proof_800}
\langle \ccalE(\bbv_2) -\ccalE(\bbv_2) , \bbv_2 -\bbv_1\rangle
  \geq  \left(1-(1-\lambda_{min}(\bbW))^{1/2}  -\alpha \rho \right) \| \bbv_2-\bbv_1\|^2.
\end{align}
Hence, if we choose the stepsize $\alpha$ such that 
\begin{equation}
\alpha \ <  \ \frac{1-(1-\lambda_{min}(\bbW))^{1/2}}{\rho} ,
\end{equation}
then $\ccalE$ is strongly monotone. In particular, if we set $\alpha \leq ({1-(1-\lambda_{min}(\bbW))^{1/2}})/({2\rho})$, then $\ccalE$ is $({1-(1-\lambda_{min}(\bbW))^{1/2}})/(2)$-strongly monotone. As $1-(1-x)^{1/2}$ is lower bounded by $x/2$ for $x\in [0,1]$, this result shows that the operator  $\ccalE=\bbD+\ccalT$ is $(\lambda_{min}(\bbW)/4)$-strongly monotone when $\alpha \leq({1-(1-\lambda_{min}(\bbW))^{1/2}})/({2\rho})$.

\subsection{Proof of Theorem~\ref{thm:approximate_cnvg}}\label{app:thm:approximate_cnvg}

	Consider the operator $\FM \defi \BM + \RM$. Note that using the update rules in (\ref{eqn_derivation_a}) and (\ref{eqn_derivation_b}), we can show that 
	\begin{align}\label{proof_thm_100}
	\alpha\FM(\zB^{t+1}) 
	&= (2\hbW-\bbI) \bbz^t  -\bbz^{t+1}- \bbU\bbq^t	\nonumber\\
	&= (2\hbW-\bbI) \bbz^t  -\bbz^{t+1}- \bbU\bbq^{t+1}+	(\bbI-\hbW) \bbz^t\nonumber\\
	&= \hbW \bbz^t  -\bbz^{t+1}- \bbU\bbq^{t+1}\nonumber\\
	&= \hbW \bbz^t  -\bbz^{t+1}- \bbU\bbq^{t+2}+	(\bbI-\hbW) \bbz^{t+1}\nonumber\\
	&= \hat{\WB}(\zB^t - \zB^{t+1}) - \bbU\qB^{t+2}.
	\end{align}
	Based on the first optimality condition in (\ref{eqn_optimality}) we know that 
	$\alpha\FM(\zB^*) + \UB\qB^* = \zeroB$. 
	By subtracting the optimality condition from \eqref{proof_thm_100}, we obtain that
	\begin{equation}\label{proof_thm_200}
		\alpha[\FM(\zB^{t+1}) - \FM(\zB^*)] = \hat{\WB}(\zB^{t} - \zB^{t+1}) - \UB(\qB^{t+2} - \qB^*).
	\end{equation}
	Using the result in \eqref{proof_thm_200}, we can write
	\begin{align}\label{444}
	&\langle \zB^{t+1} - \zB^*, \alpha [\FM(\zB^*) - \FM(\zB^{t+1})]\rangle \nonumber\\
	&\qquad = \langle \zB^{t+1} - \zB^*, -\hat{\WB}(\zB^t - \zB^{t+1}) + \UB(\qB^{t+2} - \qB^*)\rangle \nonumber\\
	&\qquad =  \langle \zB^{t+1} - \zB^*, {\hat{\WB}}(\zB^{t+1} - \zB^t)\rangle + \langle \zB^{t+1} - \zB^*, \UB(\qB^{t+2} - \qB^*)\rangle \nonumber\\
	&\qquad =  \langle \zB^{t+1} - \zB^*, {\hat{\WB}}(\zB^{t+1} - \zB^t)\rangle + \langle \qB^{t+2} - \qB^{t+1}, \qB^{t+2} - \qB^*\rangle,
	\end{align}
	where the last equality uses the definition of $\qB^t$ and that $\UB\zB^*=\zeroB$. By applying the generalized Law of cosines $2\langle a, b\rangle = \|a\|^2 + \|b\|^2 - \|a -b\|^2$, we can write the first inner product as
\begin{align}\label{0000}
 \langle \zB^{t+1} - \zB^*, {\hat{\WB}}(\zB^{t+1} - \zB^t)\rangle
 =\frac{1}{2}\left( \|\zB^{t+1} - \zB^*\|_{\hat{\WB}}^2 
 + \|\zB^{t+1} - \zB^{t}\|_{\hat{\WB}}^2
- \|\zB^{t} - \zB^*\|_{\hat{\WB}}^2\right),
\end{align}	
and the second inner product as
\begin{align}\label{1111}
\langle \qB^{t+2} - \qB^{t+1}, \qB^{t+2} - \qB^*\rangle
 =\frac{1}{2}\left( \|\qB^{t+2} - \qB^*\|^2 
 + \|\qB^{t+2} - \qB^{t+1}\|^2
- \|\qB^{t+1} - \qB^*\|^2\right),
\end{align}	
Substitute the expressions in \eqref{0000} and \eqref{1111} into \eqref{444} to obtain
	\begin{align}\label{2323}
	&\langle \zB^{t+1} - \zB^*, \alpha [\FM(\zB^*) - \FM(\zB^{t+1})]\rangle \nonumber\\
	&\quad =\frac{1}{2}\left( \|\zB^{t+1} - \zB^*\|_{\hat{\WB}}^2 
 + \|\zB^{t+1} - \zB^{t}\|_{\hat{\WB}}^2
- \|\zB^{t} - \zB^*\|_{\hat{\WB}}^2\right)\nonumber\\
&\qquad + \frac{1}{2}\left( \|\qB^{t+2} - \qB^*\|^2 
 + \|\qB^{t+2} - \qB^{t+1}\|^2
- \|\qB^{t+1} - \qB^*\|^2\right)
	\end{align}

	Now we define matrix $\bbM \in \reals^{2Nd\times 2Nd}$ and sequence of vectors $\bbphi^t\in \reals^{2Nd}$ as
	\begin{equation}\label{eqn_def_Q_X}
	 \MB \triangleq \begin{bmatrix}
	\hat{\WB} &  \ZEROBB \\
	\ZEROBB & \IB
	\end{bmatrix}
	, \ \ 
	\bbphi^t \triangleq \begin{bmatrix}
	\zB^t \\
	\qB^{t+1}
	\end{bmatrix}.
	\end{equation}
Considering these definitions we can write the inequality in \eqref{2323} as
	\begin{equation}\label{kam}
	\|\bbphi^{t+1} - \bbphi^*\|^2_{\MB} + \|\bbphi^{t+1} - \bbphi^t\|^2_{\MB} - \|\bbphi^{t} - \bbphi^*\|^2_{\MB} = 2\alpha\langle \zB^{t+1} - \zB^*, \FM(\zB^*) - \FM(\zB^{t+1})\rangle.
	\end{equation}
	Using the fact that $\ccalF=\BM + \RM$ is $\rho$-weakly monotone, we have
\begin{align}\label{kachal}
2\alpha\langle \zB^{t+1} - \zB^*, \FM(\zB^*) - \FM(\zB^{t+1})\rangle
\leq 2\alpha \rho\| \zB^{t+1} - \zB^*\|^2.
\end{align}	
	Replace the upper bound in \eqref{kachal} into \eqref{kam} and regroup the terms to obtain
	\begin{equation}
		\|\bbphi^{t} - \bbphi^*\|^2_{\MB} - \|\bbphi^{t+1} - \bbphi^*\|^2_{\MB} \geq \|\bbphi^{t+1} - \bbphi^t\|^2_{\MB} - 2\alpha \rho \|\zB^{t+1} - \zB^*\|^2.
	\end{equation}
	Replace $\|\zB^{t+1} - \zB^*\|^2 = \sum_{n=1}^{N} \|\zB^{t+1}_n - \zB^*_n\|^2 $ by its upper bound $N D^2$ where $D$ is the diameter of the set $\mathcal{Z}$.
		\begin{equation}
		\|\bbphi^{t} - \bbphi^*\|^2_{\MB} - \|\bbphi^{t+1} - \bbphi^*\|^2_{\MB} \geq \|\bbphi^{t+1} - \bbphi^t\|^2_{\MB} - 2\alpha \rho N D^2.
	\end{equation}
	Sum the above inequality from $t = 0$ to $T-1$ to obtain
	\begin{equation}
		\|\bbphi^0 - \bbphi^*\|^2_\MB - \|\bbphi^{T} - \bbphi^*\|^2_\MB \geq \sum_{t=0}^{T-1}\|\bbphi^{t+1} - \bbphi^t\|^2_{\MB} - 2\alpha \rho T N D^2
	\end{equation}
	which implies that 
		\begin{equation}
		\frac{1}{T}\sum_{t=0}^{T-1}\|\bbphi^{t+1} - \bbphi^t\|^2_{\MB}\leq  \frac{\|\bbphi^0 - \bbphi^*\|^2_\MB}{T}  + 2\alpha \rho N D^2
	\end{equation}
	Therefore, 
	\begin{align}
	\frac{1}{T}\sum_{t=0}^{T-1}\|\bbz^{t+1} - \bbz^t\|^2\leq  \frac{1}{\lambda_{\min}(\hbW)}\left(\frac{\|\bbphi^0 - \bbphi^*\|^2_\MB}{T}  + 2\alpha \rho N D^2\right)
	\end{align}
	and 
	 	\begin{align}
	\frac{1}{T}\sum_{t=0}^{T-1}\|\bbq^{t+2} - \bbq^{t+1}\|^2=\frac{1}{T}\sum_{t=0}^{T-1}\|\bbU\bbz^{t+1}\|^2\leq \frac{\|\bbphi^0 - \bbphi^*\|^2_\MB}{T}  + 2\alpha \rho N D^2
	\end{align}
	Assume that we choose one of the iterates $s$ uniformly at random from $0,\dots,T-1$, then
	\begin{equation}
	\mathbb{E}_{s}[ \|\bbz^{s+1} - \bbz^s\|^2] \leq \frac{1}{\lambda_{\min}(\hbW)}\left(\frac{\|\bbphi^0 - \bbphi^*\|^2_\MB}{T}  + 2\alpha \rho N D^2\right)
	\end{equation}
	and 
		 	\begin{align}
	\mathbb{E}_{s} [ \|\bbU\bbz^{s+1}\|^2]\leq \frac{\|\bbphi^0 - \bbphi^*\|^2_\MB}{T}  + 2\alpha \rho N D^2
	\end{align}
	Now, using the fact that $\mathbb{E} [X] \leq \sqrt{\mathbb{E} [X^2]}$ and $\sqrt{a+b} \leq \sqrt{a} + \sqrt{b}$ for positive $a,b$, we have
	\begin{equation}
	\mathbb{E}_{s}[ \|\bbz^{s+1} - \bbz^s\|] \leq \sqrt{\frac{1}{\lambda_{\min}(\hbW)}}\left(\frac{\|\bbphi^0 - \bbphi^*\|_\MB}{\sqrt{T}}  + \sqrt{2\alpha \rho N} D\right)
	\end{equation}
	and 
		 	\begin{align}
	\mathbb{E}_{s} [ \|\bbU\bbz^{s+1}\|]\leq \frac{\|\bbphi^0 - \bbphi^*\|_\MB}{\sqrt{T}}  + \sqrt{2\alpha \rho N} D
	\end{align}
	
	Now based on the expression in \eqref{proof_thm_100} we have
		\begin{align}\label{proof_thm_1000}
	\alpha\FM(\zB^{t+1}) 
	= \hat{\WB}(\zB^t - \zB^{t+1}) - \bbU\qB^{t+2}.
	\end{align}
	Multiply both sides by the the vector $ \bbone_N \otimes \bbI_d $ we obtain
		\begin{align}\label{proof_thm_10000}
	\alpha( \bbone_N \otimes \bbI_d)^\top \FM(\zB^{t+1}) 
	=( \bbone_N \otimes \bbI_d)^\top(\zB^t - \zB^{t+1}) .
	\end{align}
	which implies that 
		\begin{align}\label{proof_thm_20000}
	\alpha\sum_{n=1}^N (\ccalB_n(\bbz_n^{t+1})+\ccalR_n(\bbz_n^{t+1}))
	=\sum_{n=1}^N(\bbz_n^t - \bbz_n^{t+1}) .
	\end{align}
	Therefore, 		
	\begin{align}\label{proof_thm_30000}
	\alpha\left \|\sum_{n=1}^N\ccalB_n(\bbz_n^{t+1})+\ccalR_n(\bbz_n^{t+1}) \right\|
	&\leq \sum_{n=1}^N\|\bbz_n^t - \bbz_n^{t+1}\| \nonumber\\
	&\leq \sqrt{N \sum_{n=1}^N\|\bbz_n^t - \bbz_n^{t+1}\|^2_2}\nonumber\\
	&\leq \sqrt{N} \|\bbz^t - \bbz^{t+1}\|_2
	\end{align}
	Hence, we can show that 		
	\begin{align}\label{proof_thm_40000}
	\mathbb{E}_{s} \left[\left \|\sum_{n=1}^N\ccalB_n(\bbz_n^{s+1})+\ccalR_n(\bbz_n^{s+1}) \right\| \right]
	\leq \frac{1}{\alpha} \sqrt{\frac{N}{\lambda_{\min}(\hbW)}} \left(\frac{\|\bbphi^0 - \bbphi^*\|_\MB}{\sqrt{T}}  + \sqrt{2\alpha \rho N} D \right)
	\end{align}

\subsection{Proof of Theorem \ref{thm:convergence_DPPSP}}

Consider $\hbz=[\bbz^*;\dots;\bbz^*]$ where $\bbz^*$ is the point that satisfies the condition in Assumption~\ref{MVI_assumption}. Then we can show that 
\begin{equation}
\ccalB(\bbz)^\top (\bbz-\hbz^*)\geq 0,
\end{equation}
for any $\bbz=[\bbz_1;\dots;\bbz_N]\in \ccalZ^N$. Further, note that it can be verified that $\ccalR(\bbz)^\top (\bbz-\hbz^*)\geq 0$, as there always exists a subgradient which has a positive inner product with the vector $\bbz-\hbz^*$. Hence,
\begin{equation}
(\ccalB(\bbz)+\ccalR(\bbz))^\top (\bbz-\hbz^*)\geq 0,
\end{equation}
for any $\bbz\in \ccalZ^N$. 

Now we proceed to show that there exists a vactor $\hbq^*$ such that the vector $\hbv^*=[\hbz^*;\hbq^*]$ satisfies the condition
\begin{equation}\label{tot_MVI}
\ccalT(\bbv)^\top(\bbv-\hbv^*)\geq 0, 
\end{equation}
for all $\bbv=[\bbz;\bbq]$ such that $\bbz\in\ccalZ^N$ and $\bbq\in \reals^{Nd}$. To do so, consider $\hbq^*$ as a vector that belongs to the null space $\bbU$, i.e., $\bbU\hbq^*=\bb0$. Hence, for any $\bbz\in \ccalZ^N$ we have $\bbz^\top\bbU\hbq^*=0$. Therefore, we can write 
\begin{equation}
\alpha (\ccalB(\bbz)+\ccalR(\bbz))^\top (\bbz-\hbz^*) + \bbz^\top\bbU\hbq^*\geq 0,
\end{equation}
Furter, we know that $\hbz^*$ belongs to the null space of $\bbU$ and therefore $(\hbz^*)^\top\bbU\bbq=0$ for any $\bbq\in \reals^{Nd}$. Therefore, we have 
\begin{equation}
\alpha (\ccalB(\bbz)+\ccalR(\bbz))^\top (\bbz-\hbz^*) + \bbz^\top\bbU\hbq^* -\bbq^\top\bbU\hbz^*\geq 0,
\end{equation}
Now add and subtract $\bbz^\top\bbU\bbq$ to the left hand side to obtain
\begin{equation}
\alpha (\ccalB(\bbz)+\ccalR(\bbz))^\top (\bbz-\hbz^*) 
- \bbz^\top\bbU(\bbq-\hbq^*) +\bbq^\top\bbU(\bbz-\hbz^*)\geq 0,
\end{equation}
for any $\bbz\in\ccalZ^N$ and $\bbq\in \reals^{Nd}$. This expression can also be written as
\begin{align}
\begin{bmatrix}
	\alpha[\BM + \RM](\bbz) +\bbU\bbq \\
	-\bbU\bbz
	\end{bmatrix}^\top
	\begin{bmatrix}
	\bbz-\hbz^* \\
	\bbq-\hbq^*
	\end{bmatrix}
	\geq 0
\end{align}
which is equivalent to 
\begin{align}\label{yechi}
\left(
\begin{bmatrix}
	\alpha[\BM + \RM] & \bbU \\
	-\bbU &\bb0
	\end{bmatrix}
	\begin{bmatrix}
	\bbz\\
	\bbq
	\end{bmatrix}
	\right)^\top
	\begin{bmatrix}
	\bbz-\hbz^* \\
	\bbq-\hbq^*
	\end{bmatrix}
	\geq 0
\end{align}
By considering the definitions of $\ccalT$, $\bbv$, and $\hbv^*$, it can be verified that \eqref{yechi} implies \eqref{tot_MVI}.

Note that the sequence of iterates generated by our proposed method satisfy the condition
\begin{equation}\label{hassani}
\ccalT(\bbv^{t+1}) + \bbD(\bbv^{t+1}-\bbv^t)= \bb0,
\end{equation}
Consider the definition $\|\bba\|^2_\bbA:= \bba^\top\bbA\bba$. Then, it can be verified that
\begin{align}
\|\bbv^{t}-\hbv^*\|_\bbD^2
&= \|\bbv^{t+1}-\hbv^*\|_{\bbD}^2+ \|\bbv^{t+1}-\bbv^{t}\|_\bbD^2 +2 (\bbv^{t}-\bbv^{t+1})^T\bbD(\bbv^{t+1}-\hbv^*)
\nonumber\\
&=  \|\bbv^{t+1}-\hbv^*\|_\bbD^2+ \|\bbv^{t+1}-\bbv^{t}\|_\bbD^2 +2 \ccalT(\bbv^{t+1})^T(\bbv^{t+1}-\hbv^*)
\nonumber\\
&\geq \|\bbv^{t+1}-\hbv^*\|_\bbD^2+ \|\bbv^{t+1}-\bbv^{t}\|_\bbD^2
\end{align}
If we sum both sides from $k=0$ to $k=T-1$ we obtain that
\begin{align}
\sum_{t=0}^{T-1} \|\bbv^{t+1}-\bbv^{t}\|^2
& \leq \|\bbv^{0}-\hbv^*\|^2-\|\bbv_{T}-\hbv^*\|^2 
\end{align}
Therefore, if we choose one of the indices $k$ from $0$ to $T-1$ with probability $1/T$ the expected value of the random variable $\|\bbv^{t+1}-\bbv^{t}\|^2$ will be
\begin{align}
\mathbb{E}_s[{\|\bbv^{s+1}-\bbv^{s}\|^2}]
&=\frac{1}{T}\sum_{k=0}^{T-1} \|\bbv^{t+1}-\bbv^{t}\|^2\nonumber\\
&\leq \frac{\|\bbv^{0}-\hbv^*\|^2}{T}
\end{align}
 Based on the equality in \eqref{hassani} we can show that 
 \begin{align}
\begin{bmatrix}
	\alpha[\BM + \RM](\bbz^{t+1}) +\bbU\bbq^{t+1}\\
	-\bbU\bbz^{t+1}
	\end{bmatrix}=
	 \begin{bmatrix}
\bbz^t-\bbz^{t+1}+ \UB(\bbq^t-\bbq^{t+1}) \\
\bbU(\bbz^t-\bbz^{t+1})+ \bbq^t-\bbq^{t+1}
\end{bmatrix}
\end{align}
 Multiply both sides of the first equality by the vector $\bbone_N \otimes \bbI_d$ to obtain
 which implies that 
		\begin{align}\label{proof_thm_20001}
	\alpha\sum_{n=1}^N (\ccalB_n(\bbz_n^{t+1})+\ccalR_n(\bbz_n^{t+1}))
	=\sum_{n=1}^N(\bbz_n^t - \bbz_n^{t+1}) .
	\end{align}
	Therefore, 		
	\begin{align}\label{proof_thm_300011}
	\alpha\left \|\sum_{n=1}^N\ccalB_n(\bbz_n^{t+1})+\ccalR_n(\bbz_n^{t+1}) \right\|
	&\leq \sum_{n=1}^N\|\bbz_n^t - \bbz_n^{t+1}\| \nonumber\\
	&\leq \sqrt{N \sum_{n=1}^N\|\bbz_n^t - \bbz_n^{t+1}\|^2_2}\nonumber\\
	&\leq \sqrt{N} \|\bbz^t - \bbz^{t+1}\|_2
	\end{align}
	Hence, we can show that 		
	\begin{align}\label{proof_thm_30001}
	\mathbb{E}_{s} \left[\left \|\sum_{n=1}^N\ccalB_n(\bbz_n^{s+1})+\ccalR_n(\bbz_n^{s+1}) \right\|^2\right]
	\leq \frac{N}{\alpha^2} \frac{\|\bbv^{0}-\hbv^*\|^2}{T}
	\end{align}
	Further, we can show that
		 	\begin{align}\label{dovom}
	\mathbb{E}_{s} [ \|\bbU\bbz^{s}\|^2]
	&= \frac{1}{T}\sum_{t=0}^{T-1} \|\bbU\bbz^{s}\|^2\nonumber\\
	&= \frac{1}{T}\sum_{t=0}^{T-1} \|\bbq^{t+1}-\bbq^{t}\|^2\nonumber\\
	&\leq \frac{\|\bbv^{0}-\hbv^*\|^2}{T}
	\end{align}
	As we simply can choose $\hbq^*$ as $\hbq^*=0$ and the initial iterate is $\bbq^0 = 0$ we can simplify  $\|\bbv^{0}-\hbv^*\|^2$ as  $\|\bbz^{0}-\hbz^*\|^2$ in both \eqref{proof_thm_30001} and \eqref{dovom}.